\documentclass[11pt]{amsart}
\usepackage{amsfonts}
\usepackage{amsmath,amsthm} 
\usepackage{hyperref} 
\usepackage{latexsym}
\usepackage{array}
\usepackage{mathrsfs}
\usepackage{amssymb}
\usepackage{enumerate}
\usepackage[francais,english]{babel}
\usepackage{color}

\theoremstyle{plain}  \newtheorem{theorem}{Theorem}[section]
\newtheorem{lemma}[theorem]{Lemma}
\newtheorem{proposition}[theorem]{Proposition}
\newtheorem{corollary}[theorem]{Corollary} 
\newtheorem{definition}[theorem]{Definition} \theoremstyle{remark}
\newtheorem{remark}[theorem]{Remark}

\newcommand{\lsim}{  \lesssim   }

\newcommand{\de}{  \delta   }

\newcommand{\R}{  \mathbb{R}   }

\newcommand{\eps}{\varepsilon}

\newcommand{\C}{  \mathbb{C}   }
\newcommand{\Z}{  \mathbb{Z}   }
\newcommand{\N}{  \mathbb{N}   }

\newcommand{\M}{  \mathcal{M}   }

\newcommand{\Ca}{  \mathcal{C}   }

\newcommand{\E}{  \mathcal{E}   }
\newcommand{\Lc}{  \hat{\mathcal{E}}   }
\newcommand{\NF}{  \mathcal{NF}   }
\renewcommand{\H}{  \mathcal{H}   }
\renewcommand{\O}{  \mathcal{O}   }

\newcommand{\T}{  \mathbb{T}   }
\newcommand{\Tc}{  \mathcal{T}   }
\newcommand{\D}{  \mathcal{D}   }

\newcommand{\dd}{  \text{d}   }

\newcommand{\om}{  \omega   }

\renewcommand{\a}{  \alpha   }
\newcommand{\p}{  \partial   }
\renewcommand{\b}{  \beta   }
\renewcommand{\P}{  \mathcal{P}   }
\newcommand{\ga}{\gamma   }
\newcommand{\s}{  \sigma   }
\newcommand{\lan}{  \langle  }
\newcommand{\ran}{  \rangle  }
\newcommand{\ka}{  \kappa   }
\renewcommand{\r}{  \rho   }

\newcommand{\la}{  \lambda_a   }
\newcommand{\lb}{  \lambda_b   }

\renewcommand{\phi}{  \varphi  }
\renewcommand{\L}{  \mathcal{L}   }
\renewcommand{\S}{  \mathbb{S}   }
\newcommand{\diag}{\operatorname{diag}}
\newcommand{\meas}{\operatorname{meas}}
\newcommand{\card}{\operatorname{card}}

\newcommand{\be}{\begin{equation}}
\newcommand{\ee}{\end{equation}}
\newcommand{\ben}{\begin{equation*}}
\newcommand{\een}{\end{equation*}}
\newcommand{\ban}{\begin{align*}}
\newcommand{\ean}{\end{align*}}
\numberwithin{equation}{section}

\def\norma#1{\left\| #1\right\|}

\newcommand{\msb}{ \mathcal M_{s,\beta} }

 \author{ Beno\^it Gr\'ebert}
\address{Laboratoire de Math\'ematiques Jean Leray, Universit\'e de Nantes, UMR CNRS 6629\\
2, rue de la Houssini\`ere \\
44322 Nantes Cedex 03, France}
\email{benoit.grebert@univ-nantes.fr}
\author{Eric Paturel}
\address{Laboratoire de Math\'ematiques Jean Leray, Universit\'e de Nantes, UMR CNRS 6629\\
2, rue de la Houssini\`ere \\
44322 Nantes Cedex 03, France}
\email{eric.paturel@univ-nantes.fr}

\title[On reducibility of quantum harmonic oscillator on $\R^d$]{On  reducibility of  quantum  harmonic oscillator on $\R^d$ with  quasiperiodic in time potential.}
\begin{document}

\begin{abstract}

We prove that a linear $d$-dimensional Schr\"odinger equation on $\R^d$ with harmonic potential $|x|^2$ and small $t$-quasiperiodic potential 
$$i\partial_t u-\Delta u +|x|^2u +\eps V(t\omega,x)u=0, \quad x\in\R^d$$
reduces to an autonomous system for most values of the frequency vector $\omega\in\R^n$. As a consequence any solution of such a  linear PDE is almost periodic in time and remains bounded in all Sobolev norms.

% \begin{center}{\bf \large 18 mars 2016} \end{center}
  \end{abstract}
  
%\begin{altabstract}
% \end{altabstract}    

\subjclass{ }
\thanks{
}

\maketitle
\tableofcontents
\section{Introduction.}
\medskip

We consider the following linear Schr\"odinger equation in $\R^d$
\be\label{harmo}
i\, u_t(t,x)+(-\Delta +|x|^2)u(t,x) + \eps V(\om t,x) u(t,x) =0, \quad t\in\R,\ x\in\R^d\,.
\ee
Here $\epsilon >0$ is a small parameter and the frequency vector $\omega$ of forced oscillations is regarded as a parameter in $\mathcal D$ an open bounded subset of $\R^n$ . The function $V$ is a real multiplicative potential, which is quasiperiodic in time : namely $V$ is a continuous function of $(\phi,x)\in \T^n\times\R^d$ and $V$ is  $\H^s$ (see \eqref{Hs})  with $s>d/2$ with respect to the space variable $x\in\R^d$ and real analytic with respect to the angle variable $\phi\in\T^d$. 

%Precisely,  it analytically extends in $\phi$ to the domain $$\T_\s=\{z\in\C^n/2\pi\Z^n\mid |\Im z|<\s\}$$
%for some $\s>0$.
   We consider the previous equation as a linear non-autonomous equation in the complex Hilbert space $L^2(\R^d)$ and we prove (see Theorem \ref{thm:main} below) that it reduces to an autonomous system  for most values of the frequency vector $\om$.
   
   The general problem of reducibility for linear differential systems with time quasi periodic coefficients, 
   $\dot x =A(\om t)x$, goes back to Bogolyubov \cite{BMS69} and Moser \cite{M67}. Then there is a large literature around reducibility of finite dimensional systems by means of the KAM tools.  In particular,  the basic local result states the following : Consider the non autonomous linear system
   $$\dot x =A_0x+\eps F(\om t)x$$
   where $A_0$ and $F(\cdot)$ take values in $gl(k,\R)$, $\T^n\ni\phi\mapsto F(\phi)$ admits an analytic extension to a strip in $\C^n$ and the imaginary part of the eigenvalues of $A$ satisfy certain non resonance conditions,  then for $\eps$ small enough and for $\om$ in a Cantor set asymptotically full measure, this linear system is reducible to a constant coefficients system. This  result was then extended in many different directions (see in particular \cite{JS} \cite{E01} and \cite{K99}).
   
  Essentially our Theorem \ref{thm:main} is an infinite dimensional (i.e. $k=+\infty$) version of this basic result.
   
   Such kind of reducibility result for PDE using KAM machinery was first obtained by Bambusi \& Graffi  (see \cite{Bam-Gra}) for Schr\"odinger equation on $\R$ with a  $x^\beta$ potential, $\beta$ being strictly larger than 2.  Here we follow the more recent approach developed by  Eliasson \& Kuksin (see \cite{EK09}) for the Schr\"odinger equation on the multidimensional torus. The one dimensional case ($d=1$) was considered in \cite{GT} as a consequence of a nonlinear KAM theorem. In the present paper we extend \cite{GT}  to the multidimensional linear Schr\"odinger equation \eqref{harmo} by adapting the linear algebra tools.
   
   \medskip

We need some notations. Let   $$T=-\Delta + |x|^2=-\Delta +x_1^2+x_2^2+\cdots+x_d^2$$
 be the d-dimensional quantum harmonic oscillator. 
Its spectrum is  the sum of $d$ copies of the odd integers set, i.e. the spectrum of $T$ equals $$\hat \E:=\{d,d+2,d+4\cdots\}.$$
For $j\in\hat\E$ we denote the associated eigenspace $E_j$ whose dimension is
$$\card\ \{ (i_1, i_2,\cdots,i_d)\in(2\N-1)^d \mid i_1+i_2+\cdots+i_d=j \}:=d_j\leq j^{d-1}.$$
We denote $\{\Phi_{j,l}$,  $l=1,\cdots,d_j\}$, the basis of $E_j$ obtained by  $d$-tensor product  of Hermite functions:   $\Phi_{j,l}=\phi_{i_1} \otimes \phi_{i_2}\otimes\cdots\otimes\phi_{i_d}$ for some choice of $i_1+i_2+\cdots+i_d=j $.
Then  setting
$$\E:=\{(j,\ell)\in\hat\E\times\N\mid \ell=1,\cdots,d_j\}$$
$(\Phi_a)_{a\in\E}$ is a basis of  $L^2(\R^d)$ and denoting  
$$w_{j,\ell}=j\quad \text{ for }(j,\ell)\in\E$$
we have
$$T\Phi_a=w_a \Phi_a,\quad a\in\E.$$
We define on $\E$ an equivalence relation:
$$a\sim b\iff w_a=w_b$$
and denote by $[a]$ the equivalence class associated to $a\in\E$. We notice that
\be\label{block}\card\ [a] \leq w_a^{d-1}\,.\ee
For $s\geq 0$ an integer we define
\begin{align}\begin{split} \label{Hs}
\H^s=\{&f\in H^s(\R^d,\C) | x\mapsto {x}^{\alpha}\partial^\b f \in L^2(\R^d)\\
&\mbox{ for any } \alpha,\  \beta\in \N^d \mbox{ satisfying } 0\leq |\alpha|+|\b|\leq s \}.
\end{split}\end{align}
 We note that, for any $s\geq 0$,
$\H^s$ is the form domain of $T^{s}$ and the domain of $T^{s/2}$   (see for instance \cite{Helf84} Proposition 1.6.6) and that this allows to extend the definition of $\H^s$ to real values of $s\geq 0$. Furthermore for $s>d/2$, $\H^s$ is an algebra.\\
 To a function $u\in\H^s$ we associate the sequence $\xi$ of its Hermite coefficients by the formula $u(x)=\sum_{a\in\E}\xi_a\Phi_a(x).$ Then defining\footnote{Take care that our choose of the weight $w_a^{1/2}$ instead of $w_a$ is non standard.It is motivated by the relation \eqref{better}.  } $$\ell^2_s:=\{(\xi)_{a\in\E}\mid \sum_{a\in\E}w_a^{s}|\xi_a|^2<+\infty  \},$$ we have for $s\geq 0$
 \be\label{better}u\in\H^{s} \iff \xi\in \ell^2_s.\ee
 Then we endow both spaces with the norm
 $$\norma{u}_s=\norma{\xi}_s= (\sum_{a\in\E}w_a^{s}|\xi_a|^2)^{1/2}.$$
 If $s$ is a positive integer, we will use the fact that the norms on $\H^{s}$ are equivalently defined as $\|T^{s/2} \phi\|_{L^{2}(\R^{d})}$ and $\sum_{0\leq |\alpha|+|\b|\leq s } \|{x}^{\alpha}\partial^\b \phi \|_{L^{2}(\R^{d})}$.\\
 We finally introduce a regularity assumption on the potential $V$:
 \begin{definition}\label{admi} A potential $V: \ \T^n\times \R^d\ni(\phi, x)\mapsto V(\phi,x)\in\R$ is {\em $s$-admissible} if $\T^n\ni\phi\mapsto V(\phi,\cdot)$ is real analytic with value in $\H^s$ with
 $$\left\{ \begin{array}{lc}s \geq 0& \mbox{if}\;\; d =1\\ s > 2(d-2) & \mbox{if}\;\; d \geq 2\end{array}\right.\,.$$\end{definition}
In particular if $V$ is admissible then the map $\T^n\ni\phi\mapsto V(\phi,\cdot)\in\H^s$ analytically extends to 
 $$\T^n_\s=\{(a+ib)\in\C^n/2\pi\Z^n\mid |b|<\s  \}$$
 for some $\s>0$. Now we can state our main Theorem:
\begin{theorem}\label{thm:LS} Assume that the potential $V: \ \T^n\times \R^d\ni(\phi, x)\mapsto \R$ is {\it $s$-admissible} (see Definition \ref{admi}).
Then, there exists $\delta_{0}>0$ (depending only on $s$ and $d$) and  $\eps_*>0$ such that for all $0\leq\eps<\eps_*$ there exists $\D_{\eps} \subset [0,2\pi)^n$ satisfying 
$$\meas (\D\setminus\D_\eps)\leq \eps^{\de_0}\,,$$  such that for all 
$\omega\in \D_{\eps} $, the linear Schr\"odinger equation 
\begin{equation}\label{LS}
 i\partial_t u+(-\Delta +|x|^2)u +\eps V(t\omega,x)u=0
 \end{equation}
 reduces to a  linear equation with constant coefficients in the energy space $\H^1$.\\
 More precisely, for all $0<\de\leq \delta_0$,  there exists $\eps_0$ such that for all $0<\eps<\eps_0$ there exists $\D_{\eps} \subset [0,2\pi)^n$ satisfying 
$$\meas (\D\setminus\D_\eps)\leq \eps^\de\,,$$
and for  $\om\in\D_\eps$, there exist a linear isomorphism  $\Psi(\phi)=\Psi_{\om,\eps}(\phi)\in\L( \H^{s'})$, for $0\leq s'\leq \max(1,s)$,  unitary on $L^{2}(\R^{d})$, which analytically depends on $\phi\in\T_{\s/2}$ and a bounded Hermitian operator $W=W_{\om,\eps}\in\L(\H^{s})$ such that $t\mapsto u(t,\cdot)\in\H^{1}$ satisfies \eqref{LS} if and only if  $t\mapsto v(t,\cdot) = \Psi(\om t) u(t, \cdot)$ satisfies the autonomous equation 
 $$ i\partial_t v+(-\Delta  +|x|^2)v +\eps W(v)=0\,.$$
 Furthermore,  for all $0\leq s'\leq \max(1,s) $,
 $$\norma{\Psi(\phi)-Id}_{\L(\H^{s'},\H^{s'+2\beta})},\ \norma{\Psi(\phi)^{-1}-Id}_{\L(\H^{s'},\H^{s'+2\beta})}\leq \eps^{1-\delta/{\delta_{0}}} \quad \forall \phi\in\T^n_{\s/2}.$$ 
 On the other hand, the infinite matrix $(W_a^b )_{a,b\in\E}$ of the operator $W$ written in the Hermite basis  $(W_a^b = \int_{\R^d}\Phi_aW(\Phi_b)\dd x)$ is block diagonal, i.e.
 $$W_a^b=0 \text{ if } w_a\neq w_b$$
 and, denoting by $[V](x) = \int_{\T^d}V(\phi,x)\dd \phi$ the mean value of $V$ on the torus $\T^d$,  and by $([V]_{a}^{b})_{a,b \in \E}$  the corresponding infinite matrix, we have 
 \be\label{Wab}\norma{(W_a^b)_{a,b\in\E}-\Pi\Big(([V]_{a}^{b})_{a,b \in \E}\Big)}_{\L(\H^{s})}\leq \eps^{1/2},\ee  where $\Pi$ is the projection on the diagonal blocks.\end{theorem}

  As a consequence of our reducibility result, we prove the following corollary concerning the solutions of \eqref{harmo}. 
  \begin{corollary}\label{coro1.3}
   Assume that $(\phi,x)\mapsto V(\phi,x)$ is {\it $s$-admissible} (see Definition \ref{admi}). Let $1\leq s'\leq \max(1,s)$ and let  $u_{0}\in \H^{s'}$. Then there exists $\eps_{0}>0$ such that for all $0<\eps<\eps_{0}$ and $\om \in \D_{\eps}$,    there exists a unique solution $u \in \mathcal{C}\big(\R\,;\,\H^{s}\big)$ of \eqref{LS} such that $u(0)=u_{0}$. Moreover, $u$ is almost-periodic in time and satisfies
  \begin{equation}\label{16}
  (1-\eps C)\|u_{0}\|_{\H^{s'}}\leq \|u(t)\|_{\H^{s'}}\leq  (1+\eps C)\|u_{0}\|_{\H^{s'}}, \quad \forall \,t\in \R,
  \end{equation}
  for some $C=C(s',s,d)$.
  \end{corollary}
Another way to understand the result of Theorem \ref{thm:LS} is in term of Floquet operator (see \cite{E01} or \cite{Wang}).  Consider on $L^{2}(\T^{n}) \otimes L^2(\R^d)$ the Floquet Hamiltonian operator
 \begin{equation}\label{Floq}
 K:=i\sum_{k=1}^n\omega_k \frac{\partial}{\partial \phi_k} -\Delta  +|x|^2 +\eps V(\phi,x),
\end{equation}
then we have
 \begin{corollary}\label{coro1.4}
  Assume that $(\phi,x)\mapsto V(\phi,x)$ is {\it $s$-admissible} (see Definition \ref{admi}). There exists $\eps_{0}>0$ such that for all $0<\eps<\eps_{0}$ and $\om \in \D_{\eps}$, the spectrum of the Floquet operator $K$ is pure point.
 \end{corollary}
Let us explain our general strategy of proof of Theorem \ref{thm:LS}.\\
In the phase space $\H^{s}\times\H^{s}$ endowed with the symplectic 2-form $idu\wedge d\bar u$ equation \eqref{harmo} reads as the Hamiltonian system associated to the Hamiltonian function
\begin{align}\label{Hharmo}
H(u,\bar u)=h(u,\bar u)+\eps q(\om t, u,\bar u)
\end{align}
where
\begin{align*}
h(u,\bar u):=& \int_{\R^d} \big(|\nabla u|^2+|x|^2|u|^2\big)dx,\\
q(\om t, u,\bar u):=& \int_{\R^d} V(\om t,x)|u|^2dx.
\end{align*}
Decomposing $u$ and $\bar u$ on the basis $(\Phi_{j,l})_{(j,l)\in\E}$ of real valued functions,
$$u=\sum_{a\in\E}\xi_{a}\Phi_a,\quad \bar u=\sum_{a\in\E}\eta_{a}\Phi_a$$
the phase space $(u,\bar u)\in\H^{s}\times\H^{s}$ becomes the phase space $(\xi,\eta)\in Y_s$
$$
Y_s=\{\zeta=(\zeta_a\in\C^2,\ a\in \E)\mid \|\zeta\|_s<\infty\} 
$$
where
$$
\|\zeta\|_s^2=\sum_{a\in\E}|\zeta_a|^2 w_a^{s}.$$
%and
%$$w_{j,\ell}=j\quad \text{ for }(j,\ell)\in\E.$$
We endow $Y_s$ with the symplectic structure $id\xi\wedge d\eta$.\\
In this setting the Hamiltonians read
\begin{align*}
h=&\sum_{a\in\E}w_a \xi_a\eta_a,\\
q=& \langle \xi,Q(\om t)\eta\rangle
\end{align*}
where $Q$ is the infinite matrix whose entries are
\be\label{Q}
Q_a^b(\om t)=\int_{\R^d}V(\om t,x)\Phi_a(x)\Phi_b(x)dx
\ee
defining a linear operator on $\ell^2(\E,\C)$
and $\langle\cdot,\cdot\rangle$ is the natural pairing on $\ell^2(\E,\C)$: $\langle\xi,\eta\rangle=\sum_{a\in\E}\xi_a\eta_a$ (no complex conjugation).\\
Therefore Theorem \ref{thm:LS} is equivalent to the reducibility problem for the Hamiltonian system  associated to quadratic non autonomous Hamiltonian
\be\label{Hharmo2}\sum_{a\in\E}w_a \xi_a\eta_a+\eps  \langle \xi,Q(\om t)\eta\rangle.\ee
This reducibility is obtained by constructing a canonical change of variables close to identity such that in the new variables the Hamiltonian is autonomous and reads
$$\sum_{a\in\E}w_a \xi_a\eta_a+\eps  \langle \xi,Q_\infty\eta\rangle$$
where  $Q_\infty$ is block diagonal: $(Q_\infty)_a^b=0$ for $w_a\neq w_b$. This last condition means that, in the new variables, there is no interaction between modes of different energies, and this leads to Corollary \ref{coro1.3}. 

The proof of the reducibility theorem is based on the  following analysis already used in \cite{Bam-Gra}, \cite{EK09}, \cite{GT}: the non homogeneous Hamiltonian system
\be\label{ham-syst0}
\left\{ \begin{array}{ll}\dot \xi_a &=- i\la \xi_{a}-i\eps \left({}^{t}Q(\om t)\xi\right)_a\\ \dot \eta_a &=i\la \eta_{a}+i\eps \left(Q(\om t)\eta\right)_a\end{array}\quad a\in\E\right.\ee
is equivalent to the homogeneous system
\be\label{ham-syst-aug}
\left\{ \begin{array}{lll}\dot \xi_a &=-i\la \xi_{a}-i\eps \left({}^{t}Q(\phi)\xi\right)_a\\ \dot \eta_a &=i\la\eta_{a}+i\eps  \left(Q(\phi)\eta\right)_a\\
\dot \phi &=\om.\end{array}\quad a\in\E,\right.\ee
Consequently the canonical change of variables is constructed applying a KAM strategy to the  Hamiltonian 
$$H(y,\phi,\xi,\eta)=\om \cdot y+\sum_{a\in\E}w_a \xi_a\eta_a+\eps  \langle \xi,Q(\phi)\eta\rangle$$
in the extended phase space $\P_s=\R^n\times\T^n\times Y_s$.

\medskip

\begin{remark}\label{rem-KG}
We can also prove a similar reducibility result for the Klein Gordon equation on the sphere $\S^d$, or for the beam equation on $\T^d$,  by adapting the matrix space $\M_{s,\b}$ defined in Section 2 (see \cite{GP15}). Nevertheless, since we need a regularizing effect of the perturbation ($\b>0$ in \eqref{Msb}), in order to apply our method we cannot use it for NLS on compact domains.
\end{remark}
\begin{remark}\label{rem-KAM}
The resolution of the reducibility problem for a linear Hamiltonian PDE leads naturally to a KAM result for the corresponding nonlinear PDE. Actually the KAM procedure for nonlinear perturbations consists, roughly speaking, in linearizing  the nonlinear equation around a solution of the linear PDE and to reduce this linearized equation to a PDE with  constant coefficients. This approach is possible in the case of the Klein Gordon equation on the sphere $\S^d$ (see \cite{GP15}) or in the one dimensional case (see \cite{GT}) with analytic regularity in the space direction $x$ :  the extension to the $d$-dimensional quantum harmonic oscillator, following the realms of this paper and \cite{GP15}, is the goal of a forthcoming paper.
\end{remark}
%\begin{remark}\label{rem-KAM}
%Usually the resolution of the reducibility   problem for a linear PDE leads naturally to a KAM result for the corresponding nonlinear PDE. Actually the KAM procedure for nonlinear perturbations consists, roughly speaking, in linearizing  the nonlinear equation around a solution of the linear PDE and to reduce this linearized equation to a PDE with  constant coefficients. Nevertheless in our case, i.e. the quantum harmonic oscillator, we have to assume $V(\phi,\cdot)\in\H^{s^*}$ with $s^*>s$ to obtain the reducibility in $\H^s$. This small loss of regularity prevent to use Theorem \ref{thm:LS} in an iterative scheme. Nevertheless this approach  is possible  in the case of the Klein Gordon equation on the sphere $\S^d$ (see \cite{GP15}) or in the one dimensional case (see \cite{GT}) with analytic regularity in the space direction $x$.
%\end{remark}
\begin{remark}\label{rem-analytic}
As a difference with \cite{EK09} and \cite{GT}, we work here in spaces of {\em finite} regularity in the space variable $x$.  This allows us to get a better control of the inverse of block diagonal matrices, especially when the dimensions of the blocks are unbounded. In return, working with finite regularity in $x$ forbids any loss in this direction, at any step of the process (which is classically bypassed in the analytic case with a reduction of the analyticity strip).
\end{remark}
\medskip

\noindent {\it Acknowledgement:} The authors acknowledge the support from the projects  ANR-13-BS01-0010-03 and   ANR-15-CE40-0001-02 of the Agence Nationale de la Recherche, and Nicolas Depauw for fruitful discussions about interpolation.

\section{ Reducibility theorem.}\label{2}

In this section we state an abstract reducibility theorem for quadratic quasiperiodic in time Hamiltonians of the form
$$\sum_{a\in\E}\la \xi_a\eta_a+\eps  \langle \xi,Q(\om t)\eta\rangle.$$

\subsection{Setting}
First we need  to introduce some notations.\\
\noindent {\bf Linear space.}
Let  $s\geq 0$, we consider the complex  weighted $\ell^2$-space
$$
\ell^2_s=\{\xi=(\xi_a\in\C,\ a\in \E)\mid \|\xi\|_s<\infty\} 
$$
where $$
\|\xi\|_s^2=\sum_{a\in\L}|\xi_a|^2 w_a^{s}.$$
Then we define$$
Y_s=\ell^2_s\times\ell^2_s=\{\zeta=(\zeta_a\in\C^2,\ a\in \E)\mid \|\zeta\|_s<\infty\} 
$$
where\footnote{We provide $\C^2$ with the euclidian norm, $|\zeta_a|=|(\xi_a,\eta_a)|=\sqrt{|\xi_a|^2+|\eta_a|^2}$.}
$$
\|\zeta\|_s^2=\sum_{a\in\L}|\zeta_a|^2 w_a^{s}.$$
%In the spaces $Y_s$ acts the linear operator $J$,
%$$
%J\ :\ \{\zeta_a\}\mapsto \{\s_2\zeta_a\}, \quad \text{with } \s_2=\left(\begin{array}{cc} 0&-1\\1&0\end{array}\right).$$
 We provide the spaces $Y_s$, $s\geq 0$, with the symplectic structure $i\dd\xi\wedge\dd \eta$.  To any $C^1$-smooth function defined on a domain $\O\subset Y_s$,  corresponds the Hamiltonian equation 
 $$\left\{\begin{array}{cc}  \dot \xi &=-i\nabla_\eta f(\xi,\eta)\\ \dot \eta &=i\nabla_\xi f(\xi,\eta)\end{array}\right.
 $$
 where $\nabla f={}^t(\nabla_\xi f,\nabla_\eta f)$ is the gradient with respect to the scalar product in $Y_0$.\\
 For any $C^1$-smooth functions, $F,G$, defined on a domain $\O\subset Y_s$, we define the Poisson bracket 
$$
 \{F,G\}=i\sum_{a\in\E}\frac{\partial F}{\partial \xi_a}\frac{\partial G}{\partial \eta_a}-\frac{\partial G}{\partial \xi_a}\frac{\partial F}{\partial \eta_a}.
 $$
We will also consider the extended phase space 
 $$\P_s=\R^n\times\T^n\times Y_s\ni(y,\phi,(\xi,\eta))$$
 For any $C^1$-smooth functions, $F,G$, defined on a domain $\O\subset \P_s$, we define the extended Poisson bracket (denoted by the same symbol)
 \be\label{Poisson}
 \{F,G\}=\nabla_yF\nabla_\phi G-\nabla_yG\nabla_\phi F+i\sum_{a\in\E}\frac{\partial F}{\partial \xi_a}\frac{\partial G}{\partial \eta_a}-\frac{\partial G}{\partial \xi_a}\frac{\partial F}{\partial \eta_a}.
 \ee
 
 \smallskip
 
  \noindent {\bf Infinite matrices.}
 We denote by $\M_{s,\b}$ the set of infinite matrix $A:\L\times \L\to \C$  
% that are symmetric
%$$
%A_a^{b}=A_{b}^a,\quad \forall a,\ b\in \L$$
that satisfy
\be\label{Msb}
|A|_{s,\b} := \sup_{a,b\in \L}(w_aw_b)^\beta\left\|A_{[a]}^{[b]}\right\| \left(\frac{\sqrt{\min(w_a,w_b)}+|w_a-w_b|}{\sqrt{\min(w_a,w_b)}}\right)^{s/2}<\infty\ee
where $A_{[a]}^{[b]}$ denotes the restriction of $A$ to the block $[a]\times[b]$ and $\|\cdot\|$ denotes the operator norm. Further we denote $\M=\M_{0,0}$.
We will also need the space $\M_{s,\b}^+$ the following subspace of $\M_{s,\b}$: an infinite matrix $A\in\M$ is in $\M_{s,\b}^+$ if
$$
|A|_{s,\b+} := \sup_{a,b\in \L}\frac{(w_aw_b)^\beta}{1+|w_a-w_b|}\left\|A_{[a]}^{[b]}\right\| \left(\frac{\sqrt{\min(w_a,w_b)}+|w_a-w_b|}{\sqrt{\min(w_a,w_b)}}\right)^{s/2}<\infty$$
The following structural  lemma is proved in Appendix: 
\begin{lemma}\label{product} Let $0<\b\leq 1$ and $s\geq 0$ there exists a constant $C\equiv C(\b,s)>0$ such that
\begin{itemize}
\item[(i)]
Let $A\in \msb$ and $B\in \msb^+$. Then $AB$ and $BA$ belong to $\msb$ and
$$|AB|_{s,\b},\ |BA|_{s,\b}\leq C|A|_{s,\b}|B|_{s,\b+}.$$
\item[(ii)]
Let $A,B\in \msb^+$. Then $AB$ and $BA$ belong to $\msb^+$ and
$$|AB|_{s,\b+},\ |BA|_{s,\b+}\leq C|A|_{s,\b+}|B|_{s,\b+}.$$
\item[(iii)] Let $A\in \msb$. Then for any $t \geq 1$, $A\in\L(\ell^2_{t},\ell^2_{-t})$  and
$$\|A\xi\|_{-t}\leq C|A|_{s,\b}\|\xi\|_{t}.$$
\item[(iv)] Let $A\in \msb^+$. Then $A\in\L(\ell^2_{s'},\ell^2_{s'+2\b})$ for all $0\leq s'\leq s$  and
$$\|A\xi\|_{s'+2\b}\leq C|A|_{s,\b+}\|\xi\|_{s'}.$$
Moreover  $A \in \L(\ell^{2}_{1},\ell^{2}_{1})$ and 
$$\|A \xi \|_{1} \leq C |A|_{s,\b+}\|\xi\|_1.$$
\end{itemize}
\end{lemma}
Notice that in particular, for all $\b>0$,  matrices in $\M_{0,\b}^+$ define bounded operator on $\ell^2_1$ but, even for $s$ large, we cannot insure that $\msb\subset\L(\ell^2)$.

\medskip

\noindent{\bf Normal form:} 
\begin{definition}\label{def22}
 A matrix $Q:\ \E\times \E\to \C$ is in normal form,  and we denote $Q\in  \NF$, if
 \begin{itemize}
 \item[(i)] $Q$ is Hermitian, i.e. $Q_b^a=\overline{Q_a^b}$,
 \item[(ii)] $Q$ is block diagonal, i.e. $Q_b^a=0$ for all $w_a\neq w_b$.
 \end{itemize}
 \end{definition}
Notice that a block diagonal matrix with bounded blocks in operator norm defines a bounded operator on $\ell^2$ and thus we have $\M_{s,\b}(\D,\s)\cap \NF\subset \L(\ell^2_s)$.\\
To a matrix $Q=(Q_a^b)\in\L(\ell^2_t,\ell^2_{-t})$  we associate in a unique way a  quadratic form on $Y_s\ni (\zeta_a)_{a\in \E}=(\xi_a,\eta_a)_{a\in\E}$ by the formula
$$q(\xi,\eta)=\langle \xi, Q \eta\rangle= \sum_{a,b\in\E}\ Q_a^b\xi_a\eta_b.$$
We notice for later use that 
\be\label{poisson}
\{q_1,q_2\}(\xi,\eta)=-i\langle \xi, [Q_1,Q_2] \eta\rangle
\ee
where $$[Q_1,Q_2]=Q_1Q_2-Q_2Q_1$$ is the commutator of the two matrices $Q_1$ and $Q_2$.\\
 If $Q\in \msb$ then
\be\label{q}
\sup_{a,b\in\E}\big\|(\nabla_\xi\nabla_\eta q)_{[a]}^{[b]}\big\| \leq\frac{|Q|_{s,\b}}{(w_aw_b)^{\beta}}\left(\frac{\sqrt{\min(w_a,w_b)}}{\sqrt{\min(w_a,w_b)}+|w_a-w_b|}\right)^{s/2}.\ee

\smallskip

\noindent{\bf Parameter.}  In all the paper $\om$ will play the role of a parameter belonging to $\D_0=[0,2\pi)^n$. All the constructed functions will depend on $\om$ with $C^1$ regularity. When a function is only defined on a Cantor subset of $\D_0$ the regularity has to be understood in the Whitney sense. 

\smallskip

\noindent{\bf A class of quadratic Hamiltonians.}  
Let $s\geq 0$, $\b>0$, $\D\subset \D_0$ and $\s>0$. We denote by $\M_{s,\b}(\D,\s)$ the set of $C^1$ mappings $$\D\times\T_\s\ni(\om,\phi)\to Q(\om,\phi)\in\M_{s,\b}$$ which is real analytic in $\phi\in\T_\s:=\{\phi\in\C^n\mid |\Im \phi|<\s  \}$. This space is equipped with the norm
$$[Q]_{s,\b}^{\D,\s}=\sup_{\substack{\om\in\D,\ j=0,1\\ |\Im \phi|<\s}}  |\partial_\om^jQ(\om,\phi)|_{s,\b}\,.$$
In view of Lemma \ref{product} (iii), to a matrix $Q \in\M_{s,\b}(\D,\s)$ we can associate the quadratic form on $Y_1$ 
  $$q(\xi,\eta;\om,\phi)=\langle \xi, Q(\om,\phi) \eta\rangle\,$$
and we have 
\be\label{q2}|q(\xi,\eta;\om,\phi)|\leq [Q]_{s,\b}^{\D,\s}\norma{(\xi,\eta)}_1^2\quad \text{ for } (\xi,\eta)\in Y_1,\ \om\in\D,\ \phi\in\T_\s\,.\ee
The subspace of $\M_{s,\b}(\D,\s)$ formed by Hamiltonians $S$ such that $S(\om,\phi)\in\M_{s,\b}^+$ is denoted by $\M_{s,\b}^+(\D,\s)$ and is equipped with the norm
$$[S]_{s,\b+}^{\D,\s}=\sup_{\substack{\om\in\D,\ j=0,1\\ |\Im \phi|<\s}}  |\partial_\om^jS(\om,\phi)|_{s,\b+}\,.$$
The space of Hamiltonians $N\in\M_{s,\b}(\D,\s)$ that are independent of $\phi$ will be denoted by $\M_{s,\b}(\D)$ and is equipped with the norm
$$[N]_{s,\b}^{\D}=\sup_{{\om\in\D,\ j=0,1}}  |\partial_\om^jN(\om)|_{s,\b}\,.$$

\smallskip

\noindent{\bf Hamiltonian flow.}    To any $S\in\M_{s,\b}^+$ with $s\geq 0$ and $\b>0$ we associate the symplectic linear change of variable on $Y_s$:
$$(\xi,\eta)\mapsto (e^{-i\,{}^{t}S}\xi,e^{i S}\eta).$$
It is well defined and invertible in $\L(Y_{s'})$ for all $0\leq s'\leq \max(1,s)$ as a consequence of  Lemma \ref{product} (iv).
We note that it corresponds to the flow at time 1 generated by the quadratic Hamiltonian $(\xi,\eta)\mapsto \langle \xi, S \eta\rangle$. Notice that a necessary and sufficient condition for this flow to preserve the symmetry $\eta = \overline{\xi}$ (verified by any initial condition considered in this paper) is 
\begin{align} \label{Sherm} {}^{t}S = \overline{S} \,,
\end{align}
that is, $S$ is a hermitian matrix.\\
 When $S$ also depends smoothly on $\phi$, $\T^n\ni\phi\mapsto S(\phi)\in\M_{s,\b}^+$ we associate to $S$ the symplectic linear change of variable on the extended phase space $\P_s$: 
\be\label{flot}\Phi_S(y,\phi, \xi,\eta)\mapsto (\tilde y,\phi,e^{-i\,{}^{t}S}\xi,e^{i S}\eta)\ee
where $\tilde y$ is the solution at time $t=1$ of the equation $\dot{\tilde y}=\langle e^{-i\,{}^{t}S}\xi,\nabla_\phi Se^{i S}\eta\rangle$ with $\tilde y(0)=y.$
We note that it corresponds to the flow at time 1 generated by the  Hamiltonian $(y,\phi,\xi,\eta)\mapsto \langle \xi, S(\phi) \eta\rangle$. Concretely we will never calculate $\tilde y$ explicitly since  the non homogeneous Hamiltonian system \eqref{ham-syst0} is equivalent  to the system \eqref{ham-syst-aug} where the variable conjugated to $\phi$ is not required.

\subsection{Hypothesis on the spectrum}   

Now we formulate our hypothesis on $\la$, ${a\in\E}$:

\medskip

\noindent {\bf Hypothesis A1 -- Asymptotics.}
We assume that there exists an absolute constant $c_0>0$   such that
\be\label{laequiv}
\la\geq c_0\, w_a \quad a\in\E
\ee
and
\be\label{la-lb}
|\la-\lb|\geq {c_0}{|w_a-w_b|} \quad a,b\in\E
\ee
\noindent
{\bf Hypothesis A2 -- second Melnikov condition in measure.} There exist absolute constants $\a_1>0$, $\a_2>0$ and $C>0$ such that
the following holds:\\
for each $\ka>0$ and $K\ge1$ there exists a  closed subset $\D'=\D'(\kappa,K)\subset \D$ (where $\D$ is the initial set of vector frequencies) satisfying 
\be\label{mesomega2}
\meas(\D\setminus {\D'})\leq CK^{\a_1}  {\ka}^{\a_2}  \ee
 such that for all $\om\in{\D'}$,  all $k\in\Z^n$ with $0<|k|\leq K$ and  all $a,b\in\L$ 
 we have
 \be\label{D33}
|  k\cdot \om +\la-\lb|\geq \ka(1+|w_a-w_b|).
\ee

\subsection{The reducibility Theorem}
Let us consider  the non autonomous Hamiltonian 
\be\label{Hom}H_\om( t,\xi,\eta)=\sum_{a\in\E}\la \xi_a\eta_a+\eps  \langle \xi,Q(\om t)\eta\rangle\ee
and the associated Hamiltonian system on $ Y_s$ 
\be\label{ham-syst}
\left\{ \begin{array}{ll}\dot \xi &=- iN_0 \xi-i\eps \,{}^{t}Q(\om t)\xi\\ \dot \eta &=iN_0\eta+i\eps Q(\om t)\eta\end{array}\right.\ee
where $N_0=\diag(\la\mid a\in\E)$.
\begin{theorem}\label{thm:main} Fix $s\geq 0$, $\s>0$, $\b>0$. Assume that  $(\la)_{a\in\E}$ satisfies Hypothesis A1, A2, and that $ Q\in\M_{s,\b}(\D,\s)$. Fix $0<\de \leq\delta_{0} := \frac{\beta^{2}\alpha_{2}}{16(2+d+2\beta \alpha_{2})(d+2\beta)} $. %\leq \de_0:=\frac\b{8(d+2\b)}$.
 Then there exists $\eps_*>0$ and if $0<\eps<\eps_{*}$, there exist 
 \begin{enumerate}[(i)]
 \item a Cantor set $\D_\eps\subset \D$ with $\text{Meas}(\D\setminus \D_\eps)\leq \eps^{\de}$;
 \item a $C^1$ family
%\footnote{Here and after $C^1$ regularity  has to be understood in the Whitney sense.} 
(in $\om\in\D_\eps$) of real analytic (in $\phi\in\T_{\s/2}$) linear, unitary and symplectic coordinate transformation on $Y_0$: 
%\begin{equation*} \label{specific}
%Y_0\ni (\xi,\eta)\mapsto \Psi_\om(\phi)(\xi,\eta)\in Y_0, \quad \om\in \D_\eps,\ \phi\in \T_{\s/2},
%\end{equation*}
\begin{equation*} \label{specific}
\left\{ \begin{array}{ccl}  Y_0 &\rightarrow& Y_{0}\\  (\xi,\eta)&\mapsto& \Psi_\om(\phi)(\xi,\eta) = \langle \overline{M_{\omega}(\phi)}\xi, M_{\omega}(\phi) \eta \rangle , \quad \om\in \D_\eps,\ \phi\in \T_{\s/2}\,;\end{array}\right.
\end{equation*}
\item a $C^1$ family of quadratic autonomous Hamiltonians in normal form 
 \begin{equation*} 
\mathcal H_\om=  \langle \xi,N(\om)\eta\rangle,\quad \om\in\D_\eps\,,
 \end{equation*}
 where $N(\om)\in\NF$, in particular block diagonal (i.e. $N_a^b=0$ for $w_a\neq w_b$), and is close to $N_0=\diag(\la\mid a\in\E)$: $N(\om)-N_0\in\M_{s,\b}$ and
  \be\label{N}\| N(\om)-N_0\|_{s,\b}\leq 2\eps\quad \om\in\D_\eps\,; \ee 
 \end{enumerate}
 such that 
 \begin{equation*}
 H_\om(t,\Psi_\om(\om t)(\xi,\eta)) = \H_\om(\xi,\eta),\quad t\in\R,\ (\xi,\eta)\in Y_1,\ \om\in\D_\eps\,.
 \end{equation*}
 Furthermore $\Psi_\om(\phi)$ and $\Psi_\om(\phi)^{-1}$ are  bounded operators from $Y_{s'}$ into itself for all $0\leq s'\leq \max(1,s)$ and they are close to identity: 
\be\label{estim-L}\norma{M_\om(\phi)-Id}_{\L(\ell^{2}_{s'},\ell^{2}_{s'+2\beta})},\ \norma{M_\om(\phi)^{-1}-Id}_{\L(\ell^{2}_{s'},\ell^{2}_{s'+2\beta})}\leq \eps^{1-\de/{\delta_{0}}}\,.\ee

 %The  exponent $\a>0$  only depends  on $s$, $\b$, $n$, $d$, $c_0$, $\a_1$, $\a_2$.
% Moreover  
% $$\| N\|_{s,\b}\leq C\ ,$$
%  the new external frequencies  are close to the original ones
% \begin{equation*}
%|\mu_a-\la|\leq C\frac{\eps}{w_a^\b} ,
%\end{equation*} 
% and the new frequencies satisfy a non resonant condition, there exists $\alpha >0$ such that
%  \begin{equation*}
% \big| k\cdot \om + \mu_ar(\om) -\mu_b(\om)\big|\geq {\alpha} \  \frac{ |w_a-w_b|}{1+|k|^{\tau}},\quad k\in\Z^n,\ a,b\in\E .
% \end{equation*}
\end{theorem}
  \begin{remark}
Although $\Psi_\om(\phi)$ is defined on $ Y_0 $, the normal form $N$ (in particular $N_0$)  defines a quadratic form  on $ Y_s $ only when $s\geq 1$. Nevertheless its flow is well defined and continuous from  $Y_0$ into itself (cf. \eqref{ls3}). Fortunately our change of variable $\Psi_\om(\phi)$ is always well defined on $Y_1$ even when $ Q\in\M_{0,\b}(\D,\s)$ (i.e. when $s=0$). This is essentially a consequence of the second part of Lemma \ref{product} assertion (iv). We also remark that $\Psi_{\omega}(\phi) - Id \in \mathcal{L}(\ell^{2}_{s},\ell^{2}_{s+2\beta})$, hence it is a regularizing operator.
\end{remark}
 \begin{remark}
Notice that $\Psi_{\omega}(\phi) - Id \in \mathcal{L}(Y_{s},Y_{s+2\beta})$, i.e. it is a regularizing operator.
\end{remark}

Theorem \ref{thm:main} is proved in Section \ref{sect:kam}.

\section{Applications to the quantum harmonic oscillator on $\R^d$}
 In this section we prove Theorem \ref{thm:LS} as a corollary of Theorem \ref{thm:main}. We use notations introduced in the introduction.
 
\subsection{Verification of the hypothesis} 
 We first verify the hypothesis of Theorem \ref{thm:LS}:
 \begin{lemma}\label{OK1} When $\la=w_a$, $a\in\E$,
Hypothesis A1 and A2  hold true with $c_0 =1/2$ and $\D=[0,1]^n$.\end{lemma}
\proof
The asymptotics A1 are trivially verified with $c_0=1$.\\
It is well known (see for instance \cite{}) that
for $\tau>n$ the diophantine set
$$G_\tau(\ka):=\{\om\in[0,2\pi)^n\mid |\lan \om,k\ran+j|\geq \frac{\ka}{|k|^\tau}, \text{ for all }j\in\Z \text{ and } k\in\Z^n\setminus\{0\}\} $$
satisfies
$$\meas\big( [0,2\pi)^n\setminus G_\tau(\ka)\big)\leq C(\tau)\ka. $$
Since $w_a-w_b\in\Z$, Hypothesis A2 it satisfies choosing
$$\D=[0,1]^n, \quad \D'=G_{n+1}(\ka N^{n+1}),\quad \a_1=n+1 \text{ and }\a_2=1 .$$
\endproof

\begin{lemma}\label{OK2}
Let $d\geq 1$. Suppose that 
$$\left\{ \begin{array}{lc}s \geq 0& \mbox{if}\;\; d =1\\ s > 2(d-2) & \mbox{if}\;\; d \geq 2\end{array}\right.$$ and $V\in\H^{s}$. Then there exists  $\beta(d,s) >0$ such that  the matrix $Q$ defined by
%$s>s'$ $s> s_{0}(d)$ and $\b=\b(d,s,s')=\frac{s-s'}{4s(d+3)}$. 
$$
Q_a^b=\int_{\R^d}V(x)\Phi_a(x)\Phi_b(x)dx
$$ belongs to $\M_{s,\b(d,s)}$. Moreover, there exists $C(d,s)>0$  such that
$$|Q|_{s,\b}\leq C(d,s) \norma{V}_{s}.$$
\end{lemma}
As a consequence if $V$ is admissible (see Definition \ref{admi}) then, defining 
$$
Q_a^b(\phi)=\int_{\R^d}V(\phi,x)\Phi_a(x)\Phi_b(x)dx,
$$
the mapping $\phi\mapsto Q(\phi)$ belongs to $\msb(\D_0,\s)$ for some $\s>0$.
\proof 
First we notice that
$$\left\| Q_{[a]}^{[b]}\right\|=\sup_{\|u\|,\|v\|=1}|\langle Q_{[a]}^{[b]}u,v  \rangle|=\sup_{\substack{\Psi_a\in E_{[a]},\ \|\Psi_a\|=1\\ \Psi_b\in E_{[b]},\ \|\Psi_b\|=1}}\Big|\int_{\R^d}V(x)\Psi_a\Psi_b \, dx\Big|,$$
where $E_{[a]}$ (resp. $E_{[b]}$) is the eigenspace of $T$ associated to the cluster $[a]$ (resp. $[b]$).   
Then we follow arguments developed in \cite[Proposition 2]{Bam07} and already used in the context of the harmonic oscillator in \cite{GIP}.  The basic idea lies in the following commutator lemma: Let $A$ be a linear operator which maps $\H^s$ into itself and define the sequence of operators
$$
A_N:=[T,A_{N-1}], \quad A_0:=A$$
 then by \cite[Lemma 7]{Bam07}, we have for any $a,b\in\L$ with $w_a\neq w_b$, for any $\Psi_a\in E_{[a]}$, $\Psi_b\in E_{[b]}$ and any $N\geq 0$
$$
|{\langle A\Psi_a,\Psi_b\rangle}|\leq \frac{1}{|w_a-w_b|^N}|{\langle A_N\Psi_a,\Psi_b\rangle}| = \frac{1}{|w_a-w_b|^N} \|\Psi_{b}\|_{L^\infty}\|A_N\Psi_a \|_{L^{1}}\,.$$

Let $A$ be the  operator  given by the multiplication by the function $V(x)$.   Then, by an induction argument,
$$
A_N=\sum_{0\leq |{\alpha}|\leq N}C_{\alpha,N}D^\alpha\quad \text{with }
C_{\alpha,N}= \sum_{0\leq |{\beta}|\leq 2N-|{\alpha}|} P_{\alpha,\beta,N}(x) D^\beta V$$
and $P_{\alpha,\beta,N}$ are polynomials of degree less than $2N-|\alpha|-|\beta|$. 

We first address the case $d=1$, that we treat in the same way as in \cite{GT}. In this case, we have in \cite{KT} the following estimate on $L^{\infty}$ norm of Hermite eigenfunctions with $\|\Psi_{b}\|_{L^{2}}=1$,
$$ \|\Psi_{b}\|_{L^{\infty}} \leq w_{b}^{-1/12}\,.$$
On the other hand, for $N \geq 0$, we have
\begin{align*}
\|A_N\Psi_a \|_{L^{1}} & \leq \sum_{0\leq |{\alpha}| \leq  N} \sum_{0\leq |{\beta}| \leq 2N-|{\alpha}|} \| P_{\alpha,\beta,N}(x) D^\beta V D^\alpha \Psi_a \|_{L^1}\\
& \leq C \sum_{0\leq |{\alpha}| \leq  N} \sum_{0\leq |{\beta}| \leq 2N-|{\alpha}|}\sum_{|\gamma|\leq 2N-|\alpha|-|\beta|} \| \langle x\rangle^{\gamma} D^\beta V D^\alpha \Psi_a \|_{L^1}\\
& \leq C \sum_{0\leq |{\alpha}| \leq  N} \sum_{0\leq |{\beta}| \leq 2N-|{\alpha}|} \sum_{|\gamma|\leq  2N- |\beta|}  \| \langle x\rangle^{\gamma} D^\beta V\|_{L^{2}} \sum_{|\gamma'|\leq \alpha}\| \langle x\rangle^{-\gamma'} D^\alpha \Psi_a \|_{L^2}\\
& \leq C \|V\|_{2N} \|\Psi_{a}\|_{N}\,,
\end{align*}
where $\langle x \rangle ^{\alpha}= \Pi_{i=1}^{d}(1 +|x_{i}|^{2})^{\alpha_{i}/2}$ for $\alpha \in \N^{d}$. Moreover, since $T \Psi_a= w_a\Psi_a$ and $\|\Psi_{a}\|_{L^{2}} =1$,
\be\label{psi-s}
\norma{\Psi_a}_N\leq Cw_a^{N/2}.\ee
Therefore choosing $N=s/2$, we obtain
\begin{align}
\nonumber |{\int_{\R^d}\Psi_a\Psi_bVdx}|&\leq  \frac{C}{w_{b}^{1/12}} \Big(\frac{\sqrt{w_a}}{|w_a-w_b|}\Big)^{{s/2}} \|V\|_{s}  \\
\label{yes}&\leq \frac{2^{s/2}}{w_{b}^{1/12}}C \Big(\frac{\sqrt{w_a}}{\sqrt{w_a}+|w_a-w_b|}\Big)^{{s/2}}     \|V\|_{s}              
\end{align}
where we used that if $\sqrt{w_a}\leq |w_a-w_b|$ then $\frac{\sqrt{w_a}}{|w_a-w_b|}\leq 2 \frac{\sqrt{w_a}}{\sqrt{w_a}+|w_a-w_b|}$ while if $\sqrt{w_a}\geq |w_a-w_b|$ then $\frac{\sqrt{w_a}}{\sqrt{w_a}+|w_a-w_b|}\geq \frac 12 $ and since $|{\int_{\R^d}\Psi_a\Psi_bVdx}|\leq \norma{V}_{L^\infty}$,  \eqref{yes} is still true providing that  $C$ is large enough. Exchanging $a$ and $b$ gives
\begin{align}
|{\int_{\R^d}\Psi_a\Psi_bVdx}| &\leq \frac{2^{s/2}C}{\max(w_{a},w_{b})^{1/12}} \Big(\frac{\sqrt{\min(w_{a},w_{b})}}{\sqrt{\min(w_a,w_{b})}+|w_a-w_b|}\Big)^{{s/2}}     \|V\|_{s} \nonumber\\
& \leq \frac{2^{s/2}C}{(w_{a}w_{b})^{1/24}} \Big(\frac{\sqrt{\min(w_{a},w_{b})}}{\sqrt{\min(w_a,w_{b})}+|w_a-w_b|}\Big)^{{s/2}}     \|V\|_{s}\,, \label{interpol2}
\end{align}
hence $Q \in \M_{s, 1/24}$ and $|Q|_{s,1/24} \leq C(d,s) \|V\|_{s}$. The case $s \not \in 2 \N$ comes after a standard interpolation argument, the Stein-Weiss theorem (see e.g. \cite[Corollary 5.5.4]{BerghLofstrom}) : indeed, fixing $a$, $b$ and $s_{0} =2N$, we may estimate the norm of the linear form $V \mapsto \int_{\R^d}\Psi_a\Psi_bVdx$  acting on $\mathcal{H}^{s}$ for $s = \theta s_{0}$, $\theta \in [0,1]$, using the direct estimate
\begin{align*}
|{\int_{\R^d}\Psi_a\Psi_bVdx}| & \leq \frac{C'}{(w_{a}w_{b})^{1/24}} \|V\|_{L^{2}} 
\end{align*}
and \eqref{interpol2}, and we get
\begin{align*}
|{\int_{\R^d}\Psi_a\Psi_bVdx}| & \leq  \frac{C'}{(w_{a}w_{b})^{1/24}} \Big(\frac{\sqrt{\min(w_{a},w_{b})}}{\sqrt{\min(w_a,w_{b})}+|w_a-w_b|}\Big)^{{\theta s_{0}/2}}     \|V\|_{\theta s_{0}}\,.
\end{align*}
  \\

We now treat the case $d \geq 2$. Take $p >2$ if $d =2$ and $2 < p < \frac{2d}{d-2}$ if $d \geq 3$. Using the H\"older inequality, we get, for $\frac 1p + \frac 1q =1$,
$$|{\langle A\Psi_a,\Psi_b\rangle}|\leq \frac{1}{|w_a-w_b|^N} \|\Psi_{b}\|_{L^{p}}\|A_N\Psi_a \|_{L^{q}}\,.$$

In \cite{KT}, the $L^{p}$ estimate on Hermite eigenfunctions (with $\|\Psi_{b}\|_{L^{2}}=1$) gives
$$ \|\Psi_{b}\|_{L^{p}} \leq w_{b}^{-\tilde{\beta}(p)} \,,$$
with $\tilde{\beta}(p) = \frac{1}{3p}$ if $d=2$ (and $p \geq 10/3$) and $\tilde{\beta}(p) = \frac 12 \left(  \frac{d}{3p} - \frac{d-2}{6} \right)>0$ if $d>2$ and $ \frac{2(d+3)}{d+1} \leq p <   \frac{2d}{d-2}$. Moreover, we may estimate $\|A_N\Psi_a \|_{L^{q}}$, using Young inequality (with $\frac 12 + \frac 1r = \frac 1q$)
\begin{align*}
\|A_N\Psi_a \|_{L^{q}} & \leq  \sum_{0\leq |{\alpha}| \leq  N} \sum_{0\leq |{\beta}| \leq 2N-|{\alpha}|} \| P_{\alpha,\beta,N}(x) D^\beta V D^\alpha \Psi_a \|_{L^q}\\
& \leq C \Big( \sum_{\substack{0\leq |{\alpha}| \leq  N/2 \\ 0\leq |{\beta}| \leq 2N-|{\alpha}|}} \sum_{|\gamma|\leq 2N- \beta}   \|\langle x \rangle^{\gamma}D^{\beta}V \|_{L^{2}} \sum_{|\gamma'|\leq \alpha}\| \langle x \rangle^{-\gamma'}\Psi_a \|_{L^{r}}\\
&+ \sum_{\substack{N/2< |{\alpha}| \leq  N\\ 0\leq |{\beta}| \leq 2N-|{\alpha}|}} \sum_{|\gamma|\leq 2N - |\alpha|- |\beta|} \|\langle x \rangle^{\gamma}D^{\beta}V \|_{L^{r}} \|D^{\alpha}\Psi_a \|_{L^{2}} \Big)\\
& \leq C \Big( \|V\|_{2N} \|\Psi_{a}\|_{N/2 + \nu} + \|V\|_{3N/2+\nu} \|\Psi_{a}\|_{N} \Big)\,,
\end{align*}
using the embedding $\H^{\nu}(\R^{d}) \hookrightarrow H^{\nu}(\R^{d})$ composed with  the Sobolev embedding $H^{\nu}(\R^{d}) \hookrightarrow  L^{r}(\R^{d})$, valid for $\nu \geq d \left(\frac 12 - \frac 1r \right) = \frac dp > \frac{d-2}{2}$. Hence, for $s = 2N$ and $\nu \leq \frac N2 = \frac s 4$, i.e. $s > 2(d-2)$, we have

\begin{align*}
|{\int_{\R^d}\Psi_a\Psi_bV dx}| & \leq \frac{C_{N}}{w_{b}^{\tilde{\beta}(p)}}\frac{1}{|w_a-w_b|^{s/2}}  \norma{\Psi_a}_{s/2} \|V \|_{{s}}  \\
& \leq \frac{C_{N}}{w_{b}^{\tilde{\beta}(p)}}\frac{w_{a}^{s/4}}{|w_a-w_b|^{s/2}}  \|V \|_{{s}} \,,\end{align*}
and thus 
\begin{align}
|{\int_{\R^d}\Psi_a\Psi_bV dx}| \leq \frac{C'_{N}}{(w_{a}w_{b})^{\tilde{\beta}(p)/2}}\left(\frac{\min(w_{a},w_{b})^{1/2}}{\min(w_{a},w_{b})^{1/2} + |w_a-w_b|} \right)^{s/2} \|V \|_{{s}}\,,\label{interpol3}
\end{align}
using the same trick as in the case $d=1$. Now fixing $p(d,s)$ satisfying all the constraints $2 < p < \frac{2d}{d-2}$ and $p \geq \frac{4d}{s}$ (which is always possible since $ \frac{4d}{s} < \frac{2d}{d-2}$) and defining $\beta(d,s) = \tilde{\beta}(p(d,s))$ gives the result for an even integer $s$ satisfying $s> 2(d-2)$. In order to get the estimate for any real number $s > 2(d-2)$, we interpolate :  we take any even integer $s_{0}$ larger than $s$, and define $s_1 =0 $ and $p = +\infty$ in the case $d=2$, and $s_1=2(d-2)$, $p= \frac{2d}{d-2}$ if $d>2$. There exists $\theta \in ]0,1]$ such that $s = \theta s_{0} + (1-\theta)s_{1}$. Moreover, following the last computations, we easily find
\begin{align}
|{\int_{\R^d}\Psi_a\Psi_bV dx}| & \leq C \left(\frac{\min(w_{a},w_{b})^{1/2}}{\min(w_{a},w_{b})^{1/2} + |w_a-w_b|} \right)^{s_1/2} \|V \|_{{s_1}}\,.\label{interpol4}
\end{align}
Hence, using \cite[Corollary 5.5.4]{BerghLofstrom}, \eqref{interpol3} and \eqref{interpol4}, interpolation gives the desired estimate for $s_{1}<s \leq s_{0}$.

\endproof

\subsection{Proof of Theorem \ref{thm:LS} and Corollaries \ref{coro1.3}, \ref{coro1.4}}
The Schr\"odinger equation \eqref{LS} is a Hamiltonian system on $\H^s\times\H^s$ ($s\geq 1$) governed by the Hamiltonian function \eqref{Hharmo}. Expanding it on the orthonormal basis $(\Phi_a)_{a\in\E}$, it is    equivalent to the Hamiltonian system on $Y_s$ governed by \eqref{Hharmo2} which reads as \eqref{ham-syst} with $\la=w_a$ and $Q$ given by \eqref{Q}. By Lemmas \ref{OK1}, \ref{OK2}, if $V$ is $s$-admissible, we can apply Theorem \ref{thm:main} to \eqref{Hharmo2} and this leads to Theorem \ref{thm:LS}. More precisely, 
 in the new coordinates given by Theorem \ref{thm:LS}, $(\xi'(t),\eta'(t))=(\overline{M_{\omega}(\omega t)} \xi, M_{\omega}(\omega t) \eta)$, the system \eqref{ham-syst0} becomes autonomous and decomposes in blocks as follows (remark that since $N$ is in normal form we have ${}^{t}N = \overline{N}$):
\begin{equation}\label{ls3}
\left\{ \begin{array}{ll}     
\dot \xi'_{[a]}=-i  \overline{N}_{[a]}\xi'_{[a]} & a\in\hat\E \\
\dot \eta'_{[a]}=i N_{[a]} \eta'_{[a]}   & a\in\hat\E.
\end{array}\right.
\end{equation} 
In particular, the solution $u(t,x)$ of \eqref{LS} corresponding to the initial datum  $u_0(x)=\sum_{a\in\E} \xi(0)_a\Phi_a(x)\in\H^1$ reads $u(t,x)=\sum_{a\in\E} \xi(t)_{a}\Phi_a(x)$ with
\begin{equation}\label{fin}
 \xi(t)={}^{t} M_\om(\om t)e^{-i\overline{N}t}\,\,\overline{M_{\omega}}(0)\xi(0)\,.
 \end{equation}
  In other words, let us define the transformation $\Psi(\phi) \in \mathcal{L}(\H^{s})$ by
 $$ \Psi(\phi)\left(\sum_{a \in \E} \xi_{a} \Phi_{a}(x) \right) = \sum_{a \in \E} \left( \overline{M_{\omega}(\phi)} \xi \right)_{a} \Phi_{a}(x)\,.$$
Then $u(t,x)$ satisfies  \eqref{LS} if and only if $v(t,\cdot) = \Psi(\omega t) u(t,\cdot)$ satisfies
$$ i \partial_{t} v + (- \Delta + |x|^{2})v + \eps W(v) = 0\,,$$
where $W$ is defined as follows :
$$ W \left( \sum_{a \in \E} \xi_{a} \Phi_{a} \right) = \sum_{a \in \E} \left( N_{\omega} \xi \right)_{a} \Phi_{a}\,.$$
Furthermore, remembering the construction of $N_\om$ (see \eqref{Nom} and \eqref{tiNm}) we get that
$$\|N_\om-(N_0+\tilde N_1)\|\leq 2\eps_1=2\eps^{3/2}$$
which leads to \eqref{Wab}.
This achieves the proof of Theorem \ref{thm:LS}.

 To prove Corollary \ref{coro1.3} let us explicit the formula \eqref{fin}.
 The exponential map $e^{-i \overline{N}t}$ decomposes on the finite dimensional blocks:
 $$(e^{-i\overline N t})_{[a]}=e^{- i\overline N_{[a]}t}$$
 and $\overline N_{[a]}$  diagonalizes in orthonormal basis:
 $$P_{[a]} \overline N_{[a]}{}^tP_{[a]}= \diag(\mu_{c}),\quad P_{[a]}{}^tP_{[a]}=I_{d_a}$$ 
 where $\mu_c$ are real numbers that, in view of \eqref{N}, satisfy
 $$|\mu_a-\la|\leq C\frac{\eps}{w_a^{2\b}},\quad a\in\E.$$
 Thus 
 $$
 u(t,x)=\sum_{a\in\E} \xi_a(t)\Phi_a(x) 
 $$
 where 
 \be \label{last} \xi(t) = {}^{t}M_\om(\om t) P D(t)\, {}^{t}P\overline{M_{\omega}}(0)\xi(0) \ee
% \be\label{last}(\psi_a (t),\bar\psi_a (t))= [L_\om^{-1}(\om t)\tilde PD(t){}^t\tilde PL_\om(0)(\xi(0),\eta(0))]_a\ee
 with $$ D(t)=\diag(e^{i\mu_ct},\ c\in\E) $$ % \quad \tilde P(\xi,\eta)=(P\xi,{}^tP\eta) $$ 
 and $P$ is the $\ell^2$ unitary block diagonal map whose diagonal blocks are $P_{[a]}$. \\
In particular the solutions are all almost periodic in time with  frequencies vector $(\om,\mu)$.  Furthermore, since $\|P\xi\|_s=\|\xi\|_s$ and $M_\om(\phi)$ is close to identity (see estimate \eqref{estim-L}) we deduce \eqref{16}.\\
Now it remains to prove Corollary \ref{coro1.4}. Defining, for any $c \in \E$ the sequence $\delta^{c} \in \ell^{2}$ as $\delta^{c}_{c} = 1$ and  $\delta^{c}_{a}= 0$ if $a \neq c$, then the function $u(t,x)$ defined as
$$ u(t,x) = e^{i\mu_{c}t} \sum_{a \in [c]} \left( {}^{t}M_{\om}(\om t) P \delta^{c}\right)_{a} \Phi_{a}(x) $$
solves \eqref{LS} if and only if $\mu_{c}+k.\om$ is an eigenvalue of $K$ defined in \eqref{Floq}, with associated eigenfunction
$$ (\theta,x) \mapsto e^{i\theta.k} \sum_{a \in [c]} \left( {}^{t}M_{\om}(\theta) P \delta^{c}\right)_{a} \Phi_{a}(x)\,.$$
This shows that the spectrum of the Floquet operator \eqref{Floq} equals $\{\mu_c +k\cdot \om \mid k\in\Z^n, \ c\in\E \}$ and thus Corollary \ref{coro1.4} is proved.

% Let us  re-write \eqref{last} as
%$$\psi_a(t)=\sum_{b\in\L}\psi_{a,b}(\om t)e^{i\mu_b t}\, .$$ 
%Then we observe that $\sum_{a\in\E} \psi_{a,b}(\om t)e^{i\mu_b t} \Phi_a(x)$ solves \eqref{LS} if and only if $\mu_b+k\cdot \om$ is an eigenvalue of \eqref{Floq} (with eigenfunction  $\sum_{a\in\E} \psi_{a,b}(\phi)e^{i\phi\cdot k} \Phi_a(x)$). This shows that the spectrum of the Floquet operator \eqref{Floq} equals $\{\mu_a +k\cdot \om \mid k\in\Z^n, \ a\in\E \}$ and thus Corollary \ref{coro1.4} is proved. 

\section{Proof of Theorem \ref{thm:main}}\label{sect:kam}
\subsection{General strategy}

Let $h$ be a   Hamiltonian in normal form:
\be\label{h}
h(y,\phi,\xi,\eta)= \om\cdot y +\langle \xi, N(\om)\eta\rangle\ee
with $N$ in normal form (see Definition \ref{def22}). Notice that at the beginning of the procedure $N$ is diagonal, 
$$N=N_0=\diag(w_a,\ a\in\E)$$
 and is independent of $\om$. Let
 $q$ be a quadratic Hamiltonian of the form 
 $$q(\xi,\eta)=\langle \xi, Q(\phi)\eta\rangle$$ 
 and of size $\O (\eps)$.\\
We search for a quadratic hamiltonian  $\chi(\phi,\xi,\eta)=\langle \xi,S(\phi)\eta\rangle$ with $S=\O (\eps)$ such that its  time-one flow 
$\Phi_S\equiv \Phi_S^{t=1}$ transforms the Hamiltonian $h+ q$ into
$$
(h+ q(\phi))\circ \Phi_S=h_++ q_+(\phi),
$$
where $h_+$ is a new normal form, $\eps$-close to $h$, and the new perturbation $q_+$ is of size $\O (\eps^2)$.

As a consequence of the Hamiltonian structure we have (at least formally) that
$$(h+ q(\phi))\circ \Phi_S= h+\{ h,\chi \}+q(\phi)+ \O (\eps^2).$$
So to achieve the goal above 
we should  solve the {\it homological equation}:
\be \label{eq-homo}
\{ h,\chi \}= h_+-h -q(\phi)+\O (\eps^2).
\ee
or equivalently (see \eqref{Poisson} and \eqref{poisson})
\be\label{homo}
\om\cdot\nabla_\phi S - i[N,S]= N_+-N-Q+ \O (\eps^2).
\ee
Repeating iteratively 
the same procedure with $h_+$ instead of $h$, we will construct a change of variable $\Phi$ such that
$$
(h+ q(\phi))\circ \Phi=h_\infty\,,
$$
with $h_\infty=\om\cdot y +\langle \xi, N_\infty(\om)\eta\rangle$ in normal form.
 Note that  we will be forced to solve the homological equation, not
only for the diagonal normal form $N_0$, but for 
more general normal form Hamiltonians
\eqref{h} with  $N$ close to $N_0$ .

\subsection{Homological equation}\label{ss6.1}

 In this section we will consider a homological equation of the form
 \be\label{homo2}
 \om\cdot\nabla_\phi S-i[N,S]+Q=\text{ remainder}
 \ee  with $N$  in normal form close to $N_0$ and $Q\in\M_{s,\b}$. We will construct a solution $S\in\M^+_{s,\b}$.

\begin{proposition}\label{prop:homo}
Let $\D\subset\D_0$. Let $\D\ni\r\mapsto N(\om)\in \NF$ be a $\Ca^1$ mapping that  verifies 
\be\label{ass}
 \left\| \p_\om^j (N(\om)-N_0)_{[a]} \right\| \le \frac{c_0}{4w_a^{2\b}}\ee
for $j=0,1$, $a\in\E$  and $\om\in \D$.
Let $Q\in\M_{s,\b}$, $0<\ka\le c_0/2$ and $K\ge 1$. \\
 Then there exists a subset $\D'=\D'( \ka,K)\subset \D$, satisfying 
 \be\label{estim:D}\meas (\D\setminus \D')\leq  
 C K^{\ga_1}\ka^{\ga_2}, \ee
and there exist  $\Ca^1$-functions 
$\tilde N: \D'\to \msb \cap \NF $, $ S: \T^n_{\s}\times\D'\to \M_{s,\b}^+$ hermitian and 
 $R:\T^n_{\s}\times \D'\to \msb$,
analytic in $\phi$, such that
\be\label{homo3} 
 \om\cdot\nabla_\phi S-i[N,S]=\tilde N- Q+R
 \ee
and for all $(\phi,\om)\in \T^n_{\s'}\times \D'$, $\s'<\s$, and $j=0,1$
  \begin{align}
\label{estim-homoR}
\left| \p_\om^j R(\phi,\om)\right|_{s,\b}&\leq  C\ \frac{K^{1+\frac d2}e^{-\frac12 (\s-\s')K}}{\ka^{1+\frac{d}{2\b}} (\s-\s')^{n}}
\sup_{\substack{|\Im \phi|<\s \\ j=0,1}} |\p_\om^jQ(\phi)|_{s,\b}\,,\\
\label{estim-homoS}
 \left|\p_\om^j  S(\phi,\om)\right|_{s,\b+}&\leq C\frac{K^{d+1}}{\ka^{\frac d{\b}+2} (\s-\s')^{n}}
\sup_{\substack{|\Im \phi|<\s \\ j=0,1}} |\p_\om^jQ(\phi)|_{s,\b}\,,\\
 \label{B}
 \left|\p_\om^j \tilde N(\om)\right|_{s,\b}&\leq   \sup_{\substack{|\Im \phi|<\s \\ j=0,1}} |\p_\om^jQ(\phi)|_{s,\b}\,.\end{align} 
The constant $C$  depends on $n$, $d$, $s$, $\b$ and $|\om|$, $\ga_2=\frac{\b\a_2}{4+d+2\b\a_2}  $ and $\ga_1= \max(\a_1,2+d+n).$
\end{proposition}
\proof
Written in Fourier variables (w.r.t. $\phi$), \eqref{homo3} reads
\be\label{homo3p} 
 i\om\cdot k\  \hat S(k)-i[N,\hat S(k)]=\delta_{k,0}\tilde N- \hat Q(k)+\hat R(k)
 \ee
where  $\delta_{k,j}$ denotes the Kronecker symbol. 

 We  decompose the equation into ``components'' on each product block
$[a]\times[b]$:
\be\label{homo+}
L\, \hat S_{[a]}^{[b]}(k) 
=-i\delta_{k,0}\tilde N_{[a]}^{[b]} +i\hat Q_{[a]}^{[b]}(k){ -i}\hat R_{[a]}^{[b]}(k)
\ee
where  the operator $L:= L(k,{[a]},{[b]},\om)$ is the linear   operator, acting in the space of complex
$[a]\times[b]$-matrices defined by
$$
L\, M= \big(
 k\cdot \om \ I - N_{[a]}(\om)\big) M 
+ M N_{[b]}(\om)$$
with $N_{[a]}=N_{[a]}^{[a]}$.\\
First we solve this equation when $k=0$ and $w_a=w_b$  by defining $$\hat S_{[a]}^{[a]}(0)=0,\quad \hat R_{[a]}^{[a]}(0)=0\text{ and }\tilde N_{[a]}^{[a]}=\hat Q_{[a]}^{[a]}(0).
$$
Then we set $\tilde N_{[a]}^{[b]}=0$ for $w_a\neq w_b$ in such a way $\tilde N\in\msb\cap\NF$ and satisfies
$$
|\tilde N|_{s,\b}\leq | \hat F(0)|_{s,\b}.$$
The estimates of the derivatives with respect to $\om$ are obtained by differentiating the expressions for
$\tilde N$.

\medskip

It remains to consider the case when $k\neq 0$ or $w_a\neq w_b$.
 The matrix $N_{[a]}$ can be diagonalized in an orthonormal basis:
 $${}^tP_{[a]}N_{[a]}P_{[a]}=D_{[a]}.$$
Then we denote $\hat{S'}_{[a]}^{[b]}={}^tP_{[a]}\hat S_{[a]}^{[b]}P_{[b]}$,  $\hat{Q'}_{[a]}^{[b]}={}^tP_{[a]}\hat Q_{[a]}^{[b]}P_{[b]}$ and $\hat{R'}_{[a]}^{[b]}={}^tP_{[a]}\hat R_{[a]}^{[b]}P_{[b]}$ and we notice for later use that $\|\hat{M'}_{[a]}^{[b]}\|=\|{M}_{[a]}^{[b]}\|$ for $M=S,Q,R$.\\
In this new variables  the homological equation \eqref{homo+} reads
\be\label{homo++}(  k\cdot \om - D_{[a]})\hat{S'}_{[a]}^{[b]}(k)+{S'}_{[a]}^{[b]}(k)D_{[b]}= i\hat{Q'}_{[a]}^{[b]}(k){ -i}\hat {R'}_{[a]}^{[b]}(k).\ee
This equation can be solved term by term: let $a,b\in\E$, we set \\
\begin{align}\begin{split}\label{R'}\hat {R'}_{[a]}^{[b]}(k)&=0\quad \text{ for }|k|\leq K, \\ 
\hat{R'}_{j\ell}(k)&= \hat{Q'}_{j\ell}(k),\quad j\in[a],\ \ell\in[b],\ |k|> K,\end{split}\end{align}
and 
\begin{align}\begin{split}\label{S'}\hat{S'}_{[a]}^{[b]}(k)&=0\quad \text{ for } |k|>K \text{ or for } k=0 \text{ and } w_a=w_b,\\
\big(\hat{S'}_{[a]}^{[b]}(k)\big)_{j\ell}&=\frac {i}{ k\cdot \om \ - \alpha_j+\beta_\ell}\big(\hat{Q'}_{[a]}^{[b]}(k)\big)_{j\ell}\quad\text{in the other cases}.\end{split}\end{align}
Here $\alpha_j(\om)$ and $\beta_\ell(\om)$ denote eigenvalues of $N_{[a]}(\om)$ and $N_{[b]}(\om)$, respectively. Before the estimations of such matrices, first remark that with this resolution, we ensure that
\begin{align*}
\overline{\big(\hat{Q'}_{[a]}^{[b]}(k)\big)_{j\ell}} = \big(\hat{Q'}_{[b]}^{[a]}(-k)\big)_{\ell j} \; \Rightarrow \; \overline{\big(\hat{S'}_{[a]}^{[b]}(k)\big)_{j\ell}} = \big(\hat{S'}_{[b]}^{[a]}(-k)\big)_{\ell j}
\end{align*}
hence, if $Q'$ verifies condition \eqref{Sherm}, then this is also the case for $S'$, hence the flow induced by $S$ preserves the symmetry $\eta = \overline{\xi}$.\\ 

First notice that  \eqref{R'} classically  leads to (see for instance \cite{Kuk2})
 $$|R(\phi)|_{s,\b}=|R'(\phi)|_{s,\b}\leq C \frac{e^{-\frac12 (\s-\s')K}}{ (\s-\s')^{n}}
\sup_{|\Im\theta|<\s}| Q(\theta)|_{s,\b},\quad \text{for }|\Im\phi|<\s'.$$
In order to estimate $S$, we will use  Lemma \ref{delort} stated at the end of this section and proved in the appendix.
We face the  small divisors
\be\label{divisors}
 k\cdot\omega - \alpha_j(\om)+\beta_\ell(\om),\quad j\in [a],\ \ell\in [b].
\ee
To estimate them, we have to distinguish two cases, depending on whether $k= 0$ or not.

\medskip

{\it The case $k=0$.}
In that case, we know that $w_a\neq w_b$ 
and we use   \eqref{ass}\footnote{We use that the modulus of the eigenvalues are controlled by the operator norm of the matrix.} and \eqref{la-lb} to get 
$$|\alpha_j(\r)-\beta_\ell(\r)|\geq c_0|w_a-w_b|-\frac{c_0}{4w_a^{2\b}}
-\frac{c_0}{4w_b^{2\b}}\geq \ka(1+|w_a-w_b|).$$
This last estimate allows us to use Lemma \ref{delort} to conclude that 
\be\label{S0}
| \hat S(0)|_{\b+}\leq C \frac{1}{\ka^{1+\frac d{2\b}}}| \hat F(0)|_{\b}\,.\ee

\medskip

{\it The case $k\not=0$. }   
 Using Hypothesis A2, for any $\eta>0$,  there is a set
 $\D_1=\D( 2\eta, K)$, 
 $$\meas(\D\setminus {\D_1})\leq C K^{\a_1}\eta^{\a_2} ,$$
 such that for all $\om\in\D_1$ and  $0<|k|\leq K$ 
  $$|k\cdot\om \ -\la(\om)+\lb(\om)|\geq 2\eta(1+|w_a-w_b|).$$
 By \eqref{ass} this implies
 \begin{align*}|k\cdot\om \ -\alpha_j(\om)+\beta_\ell(\om)|&\geq 2\eta(1+|w_a-w_b|)-\frac{c_0}{4w_a^{2\b}}
 -\frac{c_0}{4w_b^{2\b}}\\
& \geq \eta(1+|w_a-w_b|)\end{align*}
 if
$$ w_b\geq w_a\geq \Big( \frac{c_0}{2\eta}\Big)^{\frac1{2\b}}.$$
Let now $ w_a\leq ( \frac{c_0}{2\eta})^{\frac1{2\b}}$. We note that $|k\cdot\om \ -\la(\om)+\lb(\om)|\leq 1$ implies that 
$w_b\leq 1+( \frac{c_0}{2\eta})^{\frac1{2\b}}+ C|k|\leq C( ( \frac{c_0}{2\eta})^{\frac1{2\b}}+ K)$.
Since $|\partial_\om(k\cdot\om)(\frac k{|k|}))|=|k|\geq 1$,  we get, using   condition \eqref{ass}, 
\be\label{estimpd}|\partial_\om(k\cdot\om -\alpha_j(\om)+\beta_\ell(\om))(\frac k{|k|}))|\geq 1/2\, .\ee
Then we  recall the following classical lemma:
\begin{lemma}\label{lem-mes}Let $f:[0,1]\mapsto\R$ a $C^1$-map satisfying $|f'(x)|\geq \delta$ for all $x\in[0,1]$ and let $\ka>0$ then
$$\meas\{x\in[0,1]\mid |f(x)|\leq \ka\}\leq \frac{\ka}{\delta}.$$
\end{lemma}
Using \eqref{estimpd} and the Lemma \ref{lem-mes}, we conclude that
\be\label{inversebis}
 |k\cdot\om -\alpha_j(\om)+\beta_\ell(\om)|\geq {\ka}(1+|w_a-w_b|)\quad \forall j\in[a],\ \forall \ell\in[b]\ee
 holds
outside a set $F_{[a],[b],k}$ of measure $\leq C w_a^dw_b^{d}(1+|w_a-w_b|)\ka$.

If $F$ is the union of  $F_{[a],[b],k}$ for $|k|\leq K$, $[a],[b]\in\hat\L$ such that $ w_a\leq  ( \frac{c_0}{2\eta})^{\frac1{2\b}}$ and $w_b \leq C( ( \frac{c_0}{2\eta})^{\frac1{2\b}}+ K)$ respectively, we have 
 \begin{align*}\meas(F)&\leq C
( \frac{c_0}{2\eta})^{\frac1{2\b}}\big(( \frac{c_0}{2\eta})^{\frac1{2\b}}+ K\big)^{d+1}K^n\big(( \frac{c_0}{2\eta})^{\frac1{2\b}}+ K\big)^{d+1}( \frac{c_0}{2\eta})^{\frac{d}{2\b}}\ka\\
&\leq C K^{n+d+2} {\eta}^{-\frac{4+d}{2\b}}\ka\,.
\end{align*}
Now we choose $\eta$ such that
$${\eta}^{\a_2}={\eta}^{-\frac{4+d}{2\b}}\ka\quad \text{i.e. }\eta=\ka^{\frac{2\b}{4+d+2\b \a_2}}.
$$
Then, as $\b\leq 1$, $\eta\geq \ka$ and we have
$$\meas(F)\leq C K^{n+d+2} \ka^{\frac{2\b\a_2}{4+d+2\b\a_2}} \,.$$
Let $\D_2=\D_1\cup F$, we have
$$\meas(\D\setminus\D_2)\leq CK^{\a_1}\eta^{\a_2}+CK^{n+d+2}\big(\frac\ka{\de_0}\big)^{\frac{2\b\a_2}{4+d+2\b\a_2}}\leq C K^{\ga_1}\ka^{\ga_2}$$
with $\ga_1=\max(\a_1,2+d+n)$, $\ga_2=\frac{2\b\a_2}{4+d+2\b\a_2}$.
Further, by construction,  for all $\r\in\D_3$,
$0<|k|\le K$, $a,b\in\L$ and $j\in[a],\ \ell\in[b]$ we have
$$|\langle k, \om(\r)\rangle \ -\alpha_j(\r)+\beta_\ell(\r)| \geq \ka(1+|w_a-w_b|).$$
 Hence using Lemma \ref{delort}  and in view of \eqref{S'}, we get  that $\hat S'(k)\in\msb^+$ and
 $$|\hat S'(k)|_{s,\b+}\leq C\frac{|\hat Q(k)|_{s,\b}K^{\frac d2}}{\ka^{1+\frac{d}{2\delta}}},\quad 0<|k|\leq K\,.$$
Combining this last estimate with \eqref{S0} we obtain a solution  $S$ satisfying for any $|\Im\phi|<\s'$
\begin{align*}
| S(\phi)|_{s,\b+}\leq &C\frac{K^{\frac d2}}{(\s-\s')^n\ka^{1+\frac{d}{2\delta}}}\sup_{|\Im \phi|<\s} |Q(\phi)|_{s,\b}\\
\end{align*}

The estimates for the derivatives with respect to $\r$ are obtained by differentiating \eqref{homo+} which leads to 
$$L(\partial_\om \hat S_{[a]}^{[b]}(k,\om))=-(\partial_\om L)\hat S_{[a]}^{[b]}(k,\om)+i\partial_\om \hat Q_{[a]}^{[b]}(k,\om)-i\partial_\om \hat R_{[a]}^{[b]}(k,\om)$$
which is an equation of the same type as \eqref{homo+} for $\partial_\om \hat S_{[a]}^{[b]}(k,\om)$ and $\partial_\om \hat R_{[a]}^{[b]}(k,\om)$ where $i\hat Q_{[a]}^{[b]}(k,\om)$ is replaced by $B_{[a]}^{[b]}(k,\om)=-(\partial_\om L)\hat S_{[a]}^{[b]}(k,\om)+i\partial_\om\hat Q_{[a]}^{[b]}(k,\om)$. This equation is solved by defining
\begin{align*}\partial_\om \hat S_{[a]}^{[b]}(k,\om)= &\chi_{|k|\le K} (k)
L(k,[a],[b],\om)^{-1}B_{[a]}^{[b]}(k,\om),\\
 \partial_\om \hat R_{[a]}^{[b]}(k,\om)=&- i{ \chi_{|k| > K} (k)}
B_{[a]}^{[b]}(k,\om)={ \chi_{|k| > K} (k)}
\partial_\rho\hat Q_{[a]}^{[b]}(k,\om)
\end{align*}
Since
$$ |(\partial_\om L)\hat S(k,\om)|_{s,\b}\leq C(K +2(\|\partial_\om A_0\|+\delta_0 )) |\hat S(k,\om)|_{s,\b}\leq CK|\hat S(k,\om)|_{s,\b}$$
we obtain
$$
|B(k,\om)|_{s,\b}\leq C K\ka^{-\frac{d}{2\b}-1}K^{d/2}|\big( |\hat Q(k)|_{s,\b}+ |\p_\om \hat Q(k)|_{s,\b}\big)
$$
and thus following the same strategy as in the resolution of \eqref{homo+} we get for $|\Im\phi|<\s'$
\begin{align*}
| \partial_\om S(\phi)|_{s,\b+}\lsim &\frac{K^{d+1}}{\ka^{\frac d{\b}+2} (\s-\s')^{n}}
\big(\sup_{|\Im \phi|<\s} |Q(\phi)|_{s,\b}+\sup_{|\Im \phi|<\s} |\p_\om Q(\phi)|_{s,\b}\big)\,,\\
|\partial_\om R(\phi)|_{s,\b}\lsim & \frac{K^{1+\frac d2}e^{-\frac12 (\s-\s')K}}{\ka^{1+\frac{d}{2\b}} (\s-\s')^{n}}
\big(\sup_{|\Im \phi|<\s} |Q(\phi)|_{s,\b}+\sup_{|\Im \phi|<\s} |\p_\om Q(\phi)|_{s,\b}\big)\,.
\end{align*}

%In this way we have constructed a solution $S_{\zeta \zeta}, R_{\zeta \zeta}, B$  of 
%the fourth component of the homological equation which satisfies all required estimates.
%To guarantee that it is real, as at the end of Section~\ref{s5.3} we replace  $S_{\zeta \zeta}, R_{\zeta \zeta}, B$
% by their real parts
%(i.e., replace $S_{\zeta \zeta}(\theta,\rho)$ by $\frac12(S_{\zeta \zeta}(\theta,\rho) + 
%\bar S_{\zeta \zeta}( \bar\theta,\rho) )$, etc.) 
 \endproof
 We end this section with the key Lemma which is an adaptation of Proposition 2.2.4 in \cite{DS1} (a similar Lemma is also proved in \cite{GP15}):
 \begin{lemma}\label{delort} Let $A\in \M$ and let $B(k)$ defined for $k\in\Z^n$ by
 \be\label{eqdelort}{B(k)}_{j}^l=\frac 1{ k\cdot \om \ - \mu_j+\mu_l}{A}_j^l,\quad j\in[a],\ \ell\in[b]\ee
 where  $\om\in\R^n$ and $(\mu_a)_{a\in\L}$ is a sequence of real numbers satisfying
 \be \label{hypdelort0}|\mu_a-w_a|\leq \min \left(\frac{C_{\mu}}{w_{a}^{\delta}}, \frac14 \right),\quad \text{ for all }a\in\L\ee
for a given $C_{\mu}>0$ and $\delta>0$,  and such that for all $a,b\in\L$ and all $|k|\leq K$
 \be\label{hypdelort} | k\cdot \om \ -\mu_a+\mu_b|\geq {\ka}(1+|w_a-w_b|).\ee
 Then $B\in\M$ and there exists a constant $C>0$ depending only on $C_{\mu}$, $|\om|$ and $\delta$  such that
 $$\|B(k)_{[a]}^{[b]}\|\leq C\frac{N^{\frac d2}}{\ka^{1+\frac{d}{2\delta}}(1+|w_a-w_b|)}\|A_{[a]}^{[b]}\|\quad \text{for all }a,b\in\L,\ |k|\leq K.$$
 \end{lemma}
The proof is based on the fact that the lemma is trivially true when $\mu_a=w_a$ is constant on each block. It is given in Appendix B.

\subsection{The KAM step.}
 Theorem \ref{thm:main} is proved by an iterative KAM procedure. We begin with the initial Hamiltonian $H_\om=h_0+q_0$ where
\be\label{h0}
h_0(y,\phi,\xi,\eta)= \om\cdot y +\langle \xi, N_0\eta\rangle\,,\ee
 $N_0=\diag (w_a,\ a\in\E)$, $\om\in\D_0$ and the quadratic perturbation $q_0(\phi,\xi,\eta)=\langle \xi,Q_0(\om,\phi)\eta\rangle$ with $Q_0=\eps Q\in\M_{s,\b}(\s_0,\D_0)$ where $\s_0=\s$. Then we construct iteratively the change of variables $\Phi_{S_m}$, the normal form $h_m=\om\cdot y +\langle \xi, N_m\eta\rangle$ and the perturbation $q_m(\phi,\xi,\eta;\om)=\langle \xi,Q_m(\om,\phi)\eta\rangle$ with $Q_m\in\M_{s,\b}(\s_m,\D_m)$ as follows: assume that the construction is done up to step $m\geq0$ then
\begin{itemize}
\item[(i)]
using Proposition \ref{prop:homo} we construct $S_{m+1}(\om,\phi)$ solution of the homological equation for $\om\in\D_{m+1}$ and $\phi\in\T^n_{\s_{m+1}}$
\be\label{homok}\om\cdot\nabla_\phi S_{m+1} - i[N_m,S_{m+1}]+Q_m= \tilde N_m+R_{m}\ee
with  $\tilde N_{m}(\om)$, $R_m(\om,\phi)$  defined for $\om\in\D_{m+1}$ and $\phi\in\T_{\s_{m+1}}$ by
\begin{align}\label{tiNm}\tilde N_{m}(\om)&=((\delta_{[j]=[\ell]}\hat Q_m(0))_{j\ell})_{j,\ell\in\E}\\
\label{Rm}R_{m}(\om,\phi)&=\sum_{|k|> K_m}\hat Q_m(\om,k)e^{ik\cdot \phi}\,;\end{align}
\item[(ii)] we define $Q_{m+1},\ N_{m+1}$ for $\om\in\D_{m+1}$ and $\phi\in\T_{\s_{m+1}}$ by
\be\label{Nm} N_{m+1}=N_m+\tilde N_m\,,\ee
and
\be\label{Qm}  Q_{m+1}=R_{m}  +\int_0^1e^{itS_{m+1}}[(1-t)(N_{m+1}-N_m+R_{m+1})+tQ_{m}, S_{m+1}]e^{-it\bar S_{m+1}}\dd t\,.
\ee
\end{itemize}
By construction, if $Q_{m}$ and $N_{m}$ are hermitian, so are $R_{m}$,  $S_{m+1}$, by the resolution of the homological equation, and  also $N_{m+1}$ and $Q_{m+1}$. 
Then we define
\begin{align}\begin{split}\label{m+1}h_{m+1}(y,\phi,\xi,\eta;\om)&=\om\cdot y +\langle \xi, N_{m+1}(\om)\eta\rangle\,, \\s_{m+1}(y,\phi,\xi,\eta;\om)&=\langle \xi, S_{m+1}(\om,\phi)\eta\rangle\,, \\
q_{m+1}(y,\phi,\xi,\eta;\om)&=\langle \xi, Q_{m+1}(\om,\phi)\eta\rangle\,.
\end{split}\end{align} 
Recall that $\Phi_{S}^t$ denotes the time $t$ flow associated to $S$ (see \eqref{flot}) and $\Phi_{S}=\Phi_{S}^1$. For any regular Hamiltonian $f$ we have, using the Taylor expansion of $g(t)=f\circ\Phi_{S_{m+1}}^t$ between $t=0$ and $t=1$
$$f\circ  \Phi_{S_{m+1}}=f+\{f,s_{m+1}\}+\int_0^1 (1-t)\{\{f,s_{m+1}\},s_{m+1}\}\circ \Phi_{S_{m+1}}^t\dd t\,.$$
Therefore we get for $\om\in\D_{m+1}$ 
\begin{align*}
(h_m+q_m)\circ  \Phi_{S_{m+1}}&=h_m+\{h_m,s_{m+1}\}+\int_0^1 (1-t)\{\{h_m,s_{m+1}\},S_{m+1}\}\circ \Phi_{S_{m+1}}^t\dd t\\
&+q_m+\int_0^1\{q_m,s_{m+1}\}\circ \Phi_{S_{m+1}}^t\dd t\\
&=h_m+\lan \xi,(\tilde N_m+R_{m})\eta\ran\\
&+\int_0^1\{(1-t)\lan \xi,(\tilde N_m+R_{m})\eta\ran+tq_m, s_{m+1}\}\circ \Phi_{S_{m+1}}^t\dd t\\
&=h_{m+1}+q_{m+1}
\end{align*}
where for the last equality we used \eqref{poisson} and \eqref{flot}.

\subsection{Iterative lemma}
Following the general scheme \eqref{homok}--\eqref{m+1} we have
$$(h_0+q_0)\circ \Phi^1_{S_{1}}\circ\cdots\circ \Phi^1_{S_m}= h_{m}+q_{m}$$
where  $q_m\in \Tc^{s,\b}(\D_m,\s_m)$  , $h_m=\om\cdot y +\lan \xi,N_m\eta\ran$ is in normal form. At step $m$  the Fourier series are truncated at order $K_m$ and the small divisors are controlled by $\ka_m$. Now we specify the choice of all the parameters for $m\geq 0$ in term of $\eps_m$ which will control with $[q_m]_{s,\b}^{\D_m,\s_m}$. \\
First we % fix
%$$ \ka_0=\eps^{1/3}.$$
% We 
define  $\eps_0=\eps$, $\s_0=\s$ and for $m\geq 1$ we choose
\begin{align*}
\s_{m-1}-\s_m=&C_* \s_0 m^{-2},\\
K_m=&2(\s_{m-1}-\s_m)^{-1}\ln \eps_m^{-1},\\
\ka_{m}=&\eps_{m-1}^{\de} 
\end{align*}
where $(C_*)^{-1} =2\sum_{j\geq 1}\frac 1{j^2}$ and $\de>0$.\\

\begin{lemma}\label{iterative} Let  $0<\de'\leq \de_0':=\frac\b{8(d+2\b)}$. There exists  $\eps_*$ depending on  $\de'$, $d$, $n$, $s$, $\b$, $\ga$, $\a_1$, $\a_2$ and $h_0$ such that, for  $0<\eps\leq\eps_*$ and
$$
\eps_{m}= \eps_0^{(3/2)^m}\quad m\geq 0\,,
$$
 we have the following:\\
For all $m\geq 1$ there exist $\D_m\subset\D_{m-1}$, $S_m\in \M_{s,\b+}(\D_m,\s_m)$, $h_m=\lan\om,y\ran +\lan \xi,N_m\eta\ran$  in normal form where $N_m\in \M_{s,\b}(\D_m)$   and there exists
$q_m\in \Tc_{s,\b}(\D_m,\s_m)$ such that for $m\geq1$
\begin{itemize}
\item[(i)]  The mapping \be \label{Phik} \Phi_{m}(\cdot,\om,\phi)=\Phi^1_{S_{m}}\ :\ Y_s \to Y_s, \quad \r\in \D_{m},\ \phi\in\T_{\s_m}\ee
is linear isomorphism  linking the Hamiltonian at step $m-1$ and the Hamiltonian at step  m, i.e.
$$(h_{m-1}+q_{m-1})\circ \Phi_{m}= h_m+q_m.$$
\item[(ii)] we have the estimates
\begin{align}
\label{DD}\meas(\D_{m-1}\setminus \D_{m})&\leq \eps_{m-1}^{\a\de'},\\
\label{NN}[\tilde N_{m-1}]_{s,\b}^{\D_m}&\leq \eps_{m-1},\\
\label{QQ}[q_m]_{s,\b}^{\D_m,\s_m}&\leq \eps_m,\\
\label{FiFi}\| \Phi_m(\cdot,\om,\phi)-Id\|_{\L(Y_s,Y_{s+2\beta})}&\leq \eps_{m-1}^{1-\nu\de'},\ \text{ for } \phi\in \T_{\s_m},\ \om\in\D_m.
\end{align}
\end{itemize}
The  exponent $\a$  and $\nu$ are given by the formulas
 $\nu= 4(\frac d\b +2)$ and $\a=\frac{\b\a_2}{2+d+2\b\a_2}  $ .
 \end{lemma}
 \proof At step 1, $h_0=\om\cdot y +\lan \xi,N_0\eta\ran$ and thus hypothesis \eqref{ass} is trivially satisfied and we can apply Proposition \ref{prop:homo} to construct $S_1$, $N_1$, $R_1$ and $\D_1$ such that for $\om\in\D_1$
 $$\om\cdot\nabla_\phi S_1-i[N_0,S_1]= N_1-N_0- Q_0+R_1.$$
Then, using \eqref{estim:D}, we have
$$\meas(\D\setminus\D_1)\leq C K_1^\ga\ka_1^{2\a}\leq  \eps_0^{\a\de'}$$
for $\eps=\eps_0$ small enough. 
Using \eqref{estim-homoS} we have for $\eps_0$ small enough
$$
[S_1]_{s,\b+}^{\D_1,\s_1}\leq C \frac{K_1^{d+1}}{\ka_1^{\frac d\b +2}(\s_0-\s_1)^n}\eps_0\leq  \eps_0^{1-\frac12\nu\de'}
$$
with $\nu= 4(\frac d\b +2)$ and thus in view of \eqref{flot} and assertion (iv) of Lemma \ref{product} we get
$$\| \Phi_1(\cdot,\om,\phi)-Id\|_{\L(Y_s,Y_{s+2\beta})}\leq  \eps_0^{1-\nu\de'}.$$
Similarly using \eqref{estim-homoR}, \eqref{B} we have
$$[N_1-N_0]_{s,\b}^{\D_1}\leq  \eps_0,$$ and
$$
[R_1]_{s,\b}^{\D_1,\s_1}\leq  \eps_0^{2-\nu\de'}
$$
for $\eps=\eps_0$ small enough. Thus using \eqref{Qm} we get
$$
[Q_1]_{s,\b}^{\D_1,\s_1}\leq C [R_1]_{s,\b}^{\D_1,\s_1}+ C([N_1-N_0]_{s,\b}^{\D_1}+[R_1]_{s,\b}^{\D_1,\s_1}+ [Q_0]_{s,\b}^{\D_1,\s_1})[S_1]_{s,\b+}^{\D_1,\s_1}\leq C\eps_0^{2-\nu\delta'}.
$$
Thus for $\de'\leq\de_0'$ and $\eps_0$ small enough
$$[Q_1]_{s,\b}^{\D_1,\s_1}\leq \eps_0^{3/2}=\eps_1.$$
 
 \medskip
 
 Now assume that we have verified Lemma \ref{iterative}  up to step $m$.
We want to perform the step $m+1$. We have $h_m=\om\cdot y +\lan \xi,N_m\eta\ran$ and since 
$$ [N_m-N_0]_{s,\b}^{\D_m}\leq [N_m-N_0]_{s,\b}^{\D_m}+\cdots+ [N_1-N_0]_{s,\b}^{\D_1}\leq \sum_{j=0}^{m-1}\eps_j\leq 2\eps_0,$$
hypothesis \eqref{ass} is satisfied and we can apply Proposition \ref{prop:homo} to construct $S_{m+1}$, $N_{m+1}$, $R_{m+1}$ and $\D_{m+1}$ such that for $\om\in\D_{m+1}$
 $$\om\cdot\nabla_\phi S_{m+1}-i[N_m,S_{m+1}]= N_{m+1}-N_m- Q_m+R_{m+1}.$$
Then, using \eqref{estim:D}, we have
$$\meas(\D_m\setminus\D_{m+1})\leq C K_{m+1}^\ga\ka_{m+1}^{2\a}\leq  \eps_m^{\a\de'}$$
for $\eps_0$ small enough. 
Using \eqref{estim-homoS} we have for $\eps_0$ small enough
$$
[S_{m+1}]_{s,\b+}^{\D_{m+1},\s_{m+1}}\leq C \frac{K_{m+1}^{d+1}}{\ka_{m+1}^{\frac d\b +2}(\s_m-\s_{m+1})^n}\eps_m\leq  \eps_m^{1-\frac12\nu\de'}.
$$
Thus in view of \eqref{flot} and assertion (iv) of Lemma \ref{product} we get
$$\| \Phi_{m+1}(\cdot,\om,\phi)-Id\|_{\L(Y_s,Y_{s+2\beta})}\leq  \eps_m^{1-\nu\de'}.$$
Similarly using \eqref{estim-homoR}, \eqref{B} we have
$$[N_{m+1}-N_m]_{s,\b}^{\D_{m+1}}\leq  \eps_m,$$ and
$$
[R_{m+1}]_{s,\b}^{\D_{m+1},\s_{m+1}}\leq  \eps_m^{2-\nu\de'}
$$
for $\eps_0$ small enough. Thus using \eqref{Qm} we get
\begin{align*}
[Q_{m+1}]_{s,\b}^{\D_{m+1},\s_{m+1}}&\leq C [R_{m+1}]_{s,\b}^{\D_{m+1},\s_{m+1}}+C\big([N_{m+1}-N_m]_{s,\b}^{\D_{m+1}}\\&+[R_{m+1}]_{s,\b}^{\D_{m+1},\s_{m+1}}+ [Q_m]_{s,\b}^{\D_{m+1},\s_{m+1}}\big)[S_{m+1}]_{s,\b+}^{\D_{m+1},\s_{m+1}}\\
&\leq C\eps_m^{2-\nu\delta'}.
\end{align*}
Thus for $\de'\leq\de_0'$ and $\eps_0$ small enough
$$[Q_{m+1}]_{s,\b}^{\D_{m+1},\s_{m+1}}\leq \eps_m^{3/2}=\eps_{m+1}.$$

\endproof

\subsection{Transition to the limit and proof of Theorem \ref{thm:main}}

Let $$\D'=\cap_{m\geq 0}\D_m.$$ In view of \eqref{DD}, this is a Borel set satisfying
$$\meas(\D\setminus\D')\leq \sum_{m\geq 0} \eps_m ^{\a\de'}\leq 2 \eps_0^{\a\de'}.$$
Let us denote $\Phi^1_N(\cdot,\om,\phi)=\Phi_{1}(\cdot,\om,\phi)\circ\cdots\circ \Phi_N(\cdot,\om,\phi)$. Due to \eqref{Phik}, it maps $Y_s$ to $Y_s$ and due to \eqref{FiFi} it satisfies for $M\leq N$ and for $\om\in\D'$, $\phi\in\T_{\s/2}$
$$\| \Phi^1_N(\cdot,\om,\phi)- \Phi^1_M(\cdot,\om,\phi)\|_{\L(Y_s,Y_{s+2\beta})}\leq  \sum_{m=M}^N\eps_m^{1-\nu\de'}\leq 2\eps_M^{1-\nu\de'}\,.$$ 
Therefore  $(\Phi_N^1(\cdot,\om,\phi))_N$ is a Cauchy sequence in $\L(Y_s,Y_{s+2\beta})$. Thus when $N\to \infty$ the maps  $\Phi^1_N(\cdot,\om,\phi)$ converge to a limit mapping $\Phi_\infty^1(\cdot,\om,\phi)\in\L(Y_s).$ Furthermore since the convergence is uniform on $\om\in\D'$ and  $\phi\in\T_{\s/2}$, $(\om,\phi)\to\Phi_\infty^1(\cdot,\om,\phi)$ is analytic in $\phi$ and $C^1$ in  $\om$. Moreover, defining $\delta = \alpha \delta'/2$ and taking $\delta_{0}=\alpha/(4\nu)$, we get
\be \label{estim-Phiinf}\|\Phi_\infty^1(\cdot,\om,\phi)-Id\|_{\L(Y_s,Y_{s+2\beta})}\leq 2\eps_0^{1-\nu\de'} < \eps_{0}^{1-\delta/{\delta_{0}}}\,.\ee
By construction, the map $\Phi_m^1(\cdot,\om,\om t)$ transforms the original Hamiltonian $$H_0=H_\om(t,\xi,\eta)=\lan \xi,N_0\eta\ran+\eps \lan\xi,Q(\om,\om t)\eta\ran$$ into $$H_m(t,\xi,\eta)=\lan \xi,N_m\eta\ran+\lan\xi,Q_m(\om,\om t)\eta\ran.$$
By \eqref{QQ}, $Q_m\to 0$ when $m\to\infty$ and by \eqref{NN}  $N_m\to N$ when $m\to\infty$ where the operator 
\begin{align}\label{Nom}
N\equiv N(\om) = N_{0} + \sum_{k=1}^{+\infty}\tilde{N}_{k}
\end{align}
is $C^1$ with respect to $\om$ and is in normal form, since this is the case for all the $N_k(\om)$.  
Further for all $\om\in\D'$ we have using \eqref{NN}
$$\norma{N(\om)-N_0}_{s,\b}\leq \sum_{m=0}^\infty \eps^m\leq 2 \eps .$$
Let us denote $\Psi_\om(\phi)=\Phi_\infty^1(\cdot,\om,\phi)$. By construction, 
$$ \Psi_{\omega}(\phi) = \langle \overline{M_{\omega}(\phi)} \xi, M_{\omega}(\phi)\eta \rangle\,,$$
where
$$ M_{\omega}(\phi) = \lim_{j \rightarrow +\infty} e^{iS_{1}(\omega,\phi)}\ldots e^{iS_{j}(\omega, \phi)}\,.$$

Further, denoting the limiting Hamiltonian $\H_\om(\xi,\eta)=\lan \xi,N\eta\ran$ we have 
$$H_\om(t,\Psi_\om(\om t)(\xi,\eta)) = \H_\om(\xi,\eta),\quad t\in\R,\ (\xi,\eta)\in Y_s,\ \om\in\D_\eps\,.$$

This concludes the proof of Theorem \ref{thm:main}.

\appendix
\section{Proof of Lemma \ref{product}}
We start with two auxiliary lemmas
\begin{lemma}\label{lem-A1}
Let $j,k,\ell\in\N\setminus\{0\}$ then
\be\label{AA}\frac{\sqrt{\min(j,k)}}{\sqrt{\min(j,k)}+|j-k|}\ \frac{\sqrt{\min(\ell,k)}}{\sqrt{\min(\ell,k)}+|\ell-k|}\leq \frac{\sqrt{\min(j,\ell)}}{\sqrt{\min(j,\ell)}+|j-\ell|}.\ee
\end{lemma}
\proof
Without loss of generality we can assume $j\leq\ell$.\\
If $k\leq j$ then $|k-\ell|\geq |j-\ell|$  and thus
\begin{align*}\frac{\sqrt{\min(j,\ell)}}{\sqrt{\min(j,\ell)}+|j-\ell|}&=\frac{\sqrt j}{\sqrt j+|j-\ell|}\geq \frac{\sqrt j}{\sqrt j+|k-\ell|}\\ &\geq \frac{\sqrt k}{\sqrt k+|k-\ell|}=\frac{\sqrt{\min(k,\ell)}}{\sqrt{\min(k,\ell)}+|k-\ell|}\end{align*}
which leads to \eqref{AA}. The case $\ell\leq k$ is similar.\\
In the case $j\leq k\leq \ell$ we have 
\begin{align*}&\frac{\sqrt{\min(j,k)}}{\sqrt{\min(j,k)}+|j-k|}\ \frac{\sqrt{\min(\ell,k)}}{\sqrt{\min(\ell,k)}+|\ell-k|}\\
 &\leq \frac{\sqrt j}{\sqrt j+|j-k|}\frac{\sqrt k}{\sqrt k+|k-\ell|}\leq \frac{\sqrt j}{\sqrt j+|j-k|}\frac{\sqrt j}{\sqrt j+|k-\ell|}\\
 &\leq\frac{\sqrt j}{\sqrt j+|j-k|+|k-\ell|}
\leq \frac{\sqrt j}{\sqrt j+|j-\ell|}=\frac{\sqrt{\min(j,\ell)}}{\sqrt{\min(j,\ell)}+|j-\ell|}\,.\end{align*}
\endproof

\begin{lemma}\label{lem-A2}
Let $j\in\N$ then
$$\sum_{k\in \N}\frac{1}{k^\b(1+|k-j|)}\leq C(\b)$$
for a constant $C(\b)>0$ depending only on $\b>0$.
\end{lemma}
\proof
We note that $$\sum_{k\in \N}\frac{1}{k^\b(1+|k-j|)}=a\star b(j)$$ where  $a_k=\frac 1k$ for $k\geq1$, $a_k=0$ for $k\leq 0$ and $b_k=\frac 1{1+|k|}$, $k\in\Z$.
We have that $b\in\ell^p$ for any $1<p\leq+\infty$ and that $a\in\ell^q$ for any $\frac1\b<q\leq+\infty$. Thus by Young inequality $a\star b\in\ell_r$  for $r$ such that $\frac 1p+\frac 1q=1+\frac 1r$. In particular choosing $q=\frac 2\b$ and $p=\frac 2{2-\b}$ we conclude that  $a\star b\in\ell_\infty$.
\endproof

\proof {\it of Lemma \ref{product}.}\\
(i) Let $a,b\in\E$
\begin{align*}
\left\| (AB)_{[a]}^{[b]} \right\|&\leq \sum_{c\in\Lc}\left\| A_{[a]}^{[c]} \right\|\left\| B_{[c]}^{[b]} \right\|\\
&\leq  \frac{|A|_{s,\b+}|B|_{s,\b}}{(w_aw_b)^{\b}}\left(\frac{\sqrt{\min(w_a,w_b)}}{\sqrt{\min(w_a,w_b)}+|w_a-w_b|}\right)^{ s/2}
\sum_{c\in\Lc}\frac 1 {w_c^{2\b}(1+|w_a-w_c|)}\\
&\leq C \frac{|A|_{\b+}|B|_{\b}}{(w_aw_b)^{\b}}\left(\frac{\sqrt{\min(w_a,w_b)}}{\sqrt{\min(w_a,w_b)}+|w_a-w_b|}\right)^{s/2}
\end{align*}
where we used that by Lemma \ref{lem-A1}
$$\frac{\sqrt{\min(w_a,w_b)}}{\sqrt{\min(w_a,w_b)}+|w_a-w_b|}\geq  \frac{\sqrt{\min(w_a,w_c)}}{\sqrt{\min(w_a,w_c)}+|w_a-w_c|}\frac{\sqrt{\min(w_c,w_b)}}{\sqrt{\min(w_c,w_b)}+|w_c-w_b|}$$
 and that by Lemma \ref{lem-A2}, $\sum_{c\in\Lc}\frac 1 {w_c^{2\b}(1+|w_a-w_c|)}\leq C$ where $C$ only depends on $\b$.\\ 
(ii) Similarly let $a,b\in\L$ and assume without loss of generality that $w_a\leq w_b$
\begin{align*}
&\left\| (AB)_{[a]}^{[b]} \right\|\leq \sum_{c\in\Lc}\left\| A_{[a]}^{[c]} \right\|\left\| B_{[c]}^{[b]} \right\|\\
&\leq \frac{|A|_{s,\b+}|B|_{s,\b+}}{(w_aw_b)^{\b}}\left(\frac{\sqrt{\min(w_a,w_b)}}{\sqrt{\min(w_a,w_b)}+|w_a-w_b|}\right)^{s/2}
\sum_{c\in\Lc}\frac 1 {w_c^{2\b}(1+|w_a-w_c|)(1+|w_b-w_c|)}\\
&\leq  \frac{2|A|_{s,\b+}|B|_{s,\b+}}{(w_aw_b)^{\b}(1+|w_a-w_b|)}\left(\frac{\sqrt{\min(w_a,w_b)}}{\sqrt{\min(w_a,w_b)}+|w_a-w_b|}\right)^{s/2}\\
&\Big(\sum_{\substack{c\in\Lc \\ w_c\leq \frac12(w_a+w_b) }}\frac 1 {w_c^{2\b}(1+|w_a-w_c|)}
 +\sum_{\substack{c\in\Lc \\  w_c\geq \frac12(w_a+w_b) }}\frac 1 {w_c^{2\b}(1+|w_b-w_c|)}\Big)\\
&\leq C\frac{|A|_{s,\b+}|B|_{s,\b+}}{(w_aw_b)^{\b}(1+|w_a-w_b|)}\left(\frac{\sqrt{\min(w_a,w_b)}}{\sqrt{\min(w_a,w_b)}+|w_a-w_b|}\right)^{s/2}\,.
\end{align*}
(iii)Let $\xi\in \ell^2_{t}$, with $t \geq1$. We have
\begin{align*}
\| A\xi\|^2_{-t}&\leq\sum_{a\in\Lc}w_a^{-t}\big(\sum_{b\in\Lc}\|A_{[a]}^{[b]}\| \|\xi_{[b]}\|\big)^2\\
&\leq |A|^2_{s,\b}\sum_{a\in\Lc}\Big(\sum_{b\in\Lc}\frac{\|w_b^{t/2}\xi_{[b]}\|}{w_a^{t/2+\b}w_b^{t/2+\b}}\left(\frac{\sqrt{\min(w_a,w_b)}}{\sqrt{\min(w_a,w_b)}+|w_a-w_b|}\right)^{s/2}\Big)^2\\
&\leq \sum_{a\in\Lc} \frac1{w_a^{t+2\b}}\sum_{b\in\Lc} \frac1{w_b^{t+2\b}}|A|^2_{s,\b}\norma{\xi}^2_{t}\, .
\end{align*}
(iv) Let $\xi\in \ell^2_s$.  We have
\begin{align*}
\| A\xi\|^2_{s+2\b}&\leq\sum_{a\in\Lc}w_a^{s+2\b}\big(\sum_{b\in\Lc}\|A_{[a]}^{[b]}\| \|\xi_{[b]}\|\big)^2\\
&\leq |A|^2_{s,\b+}\sum_{a\in\Lc}\Big(\sum_{b\in\Lc}\frac{w_a^{s/2}\|w_b^{s/2}\xi_{[b]}\|}{w_b^{s/2+\b}(1+|w_a-w_b|)}\left(\frac{\sqrt{\min(w_a,w_b)}}{\sqrt{\min(w_a,w_b)}+|w_a-w_b|}\right)^{s/2}\Big)^2\\
&\leq 2^{s+1} |A|^2_{s,\b+}\sum_{a\in\Lc}\Big(\sum_{\substack{b\in\Lc \\ w_a\leq 2w_b }}\frac{\|w_b^{s/2}\xi_{[b]}\|}{w_b^\b(1+|w_a-w_b|)}\\
&\hspace{1cm}+\sum_{\substack{b\in\Lc \\ w_a\geq 2 w_b }} \frac{\|w_b^{s/2}\xi_{[b]}\|\min(w_a,w_b)^{\frac{s}2}}{w_b^{s/2+\b}(1+|w_a-w_b|)} \Big)^2\\
&\leq 2^{s+1} |A|^2_{s,\b+}\sum_{a\in\Lc}\big(\sum_{b\in\Lc}\frac{\|w_b^{s/2}\xi_{[b]}\|}{w_b^\b(1+|w_a-w_b|)}\big)^2\,.
\end{align*}
Then we note that
$$\sum_{b\in\Lc}\frac{\|w_b^s\xi_{[b]}\|}{w_b^\b(1+|w_a-w_b|)}=u\star v(a)$$
with $u_b=\|w_b^{s/2-\beta}\xi_{[b]}\|$ and $v_b=\frac{1}{(1+|w_b|)}$. Using the Cauchy Schwarz inequality we get
$$\sum_{b\in\Lc} u_b^p\leq \big(\sum_{b\in\Lc} \|w_b^{s/2}\xi_{[b]}\|^2\big)^{/2}  \big(\sum_{b\in\Lc}w_b^{\frac{2\b p}{2-p}}  \big)^{\frac{2-p}{p}}.$$
Choosing $p=\frac2{1+\b}$ we have  $\frac{2\b p}{2-p}=2>1$ and thus  $u\in \ell^p$. Choosing  
$q=\frac{2p}{3p-2}$ we have $q=\frac{2}{2-\b}>1$ and thus  $v\in\ell^q$. Since $1/p+1/q=3/2$ we conclude that $u\star v\in\ell^2$ and
$$\norma{u\star v}_{\ell^2}\leq C\norma{u}_{\ell^p}\norma{v}_{\ell^q}.$$
This leads to the first part of (iv) since $\norma{u}_{\ell^p}\leq C \|\xi\|_s.$
Now we prove the second assertion of (iv) in a similar way : let $\xi \in \ell^{2}_{1}$, we have
\begin{align*}
\| A\xi\|^2_{1} &\leq\sum_{a\in\Lc}w_a \big(\sum_{b\in\Lc}\|A_{[a]}^{[b]}\| \|\xi_{[b]}\|\big)^2\\
&\leq |A|^2_{s,\b+}\sum_{a\in\Lc}\Big(\sum_{b\in\Lc}\frac{w_a^{1/2}\|w_b^{1/2}\xi_{[b]}\|}{(w_{a}w_{b})^{\beta} w_b^{1/2}(1+|w_a-w_b|)}\left(\frac{\sqrt{\min(w_a,w_b)}}{\sqrt{\min(w_a,w_b)}+|w_a-w_b|}\right)^{s/2}\Big)^2\\
&\leq 2^{s+1} |A|^2_{s,\b+}\sum_{a\in\Lc}\Big(\sum_{\substack{b\in\Lc \\ w_a\leq 2w_b }}\frac{\|w_b^{1/2}\xi_{[b]}\|}{(w_{a}w_b)^\b(1+|w_a-w_b|)}\\
&\hspace{1cm}+\sum_{\substack{b\in\Lc \\ w_a\geq 2 w_b }} \frac{\|w_b^{1/2}\xi_{[b]} \|w_{a}^{(1-s)/2}}{(w_{a}w_{b})^{\b} w_{b}^{1/2-s/4}(1+|w_{a}-w_{b}|)} \Big)^{2}\\
%& \leq 2^{s+2} |A|^2_{s,\b+}\sum_{a\in\Lc} \Big( \big(\sum_{b\in\Lc}\frac{\|w_b^{1/2}\xi_{[b]}\|}{w_b^\b(1+|w_a-w_b|)}\big)^2\\
%&\hspace{1cm} + \big(\sum_{\substack{b\in\Lc \\ w_a\geq 2 w_b }} \frac{\|w_b^{1/2}\xi_{[b]}}{(w_{a}w_{b})^{\b} w_{b}^{1/2-s/4}(1+|w_{a}-w_{b}|)^{1/2+s/2}} \big)^{2}\Big)\,.\\
\end{align*}
The last sum may be bounded above by (notice that $|w_a-w_b|\geq w_b$)
\begin{align*}
\sum_{\substack{b\in\Lc \\ w_a\geq 2 w_b }} \frac{\|w_b^{1/2}\xi_{[b]} \|w_{a}^{(1-s)/2}}{(w_{a}w_{b})^{\b} w_{b}^{1/2-s/4}(1+|w_{a}-w_{b}|)} & \leq \sum_{\substack{b\in\Lc \\ w_a\geq 2 w_b }} \frac{\|w_b^{1/2}\xi_{[b]}\|}{(w_{a}w_{b})^{\b} w_{b}^{1/2-s/4}(1+|w_{a}-w_{b}|)^{1/2+s/2}} \\
& \leq \frac{1}{w_{a}^{\b}}  \sum_{b\in\Lc} \frac{\|w_b^{1/2}\xi_{[b]}\|}{w_{b}^{1/2+\b/2} (1+|w_{a}-w_{b}|)^{1/2+\b/2}}\,,
\end{align*}
and this last sum is the convolution product $u'\star v' (a)$, with $u'_{b} = \frac{\|w_b^{1/2}\xi_{[b]}\|}{w_{b}^{1/2+\b}}$, which  defines a $\ell^{1}$ sequence thanks to Cauchy Schwarz inequality, and $v'_{b} = \frac{1}{(1+ w_{b})^{1/2+\b/2}}$, which defines a $\ell^{2}$ sequence.  Therefore, it is a $\ell^{2}$ sequence with index $a$. We treat the first sum in the same way as before, and we obtain
$$ \| A\xi\|^2_{1}\leq C |A|^2_{s,\b+} \|\xi\|_{1}^{2}\,.$$

\endproof

\section{Proof of Lemma \ref{delort}}

%We first remark that the claimed property only concerns the operator norms of the blocks $B_{[a]}^{[b]}$, which can be computed separately. Let $K_{1}$ and $K_{2}$ be positive integers that will be fixed later.  We define the following decomposition in $\msb$, according to the weights $w_{a}$ and $w_{b}$ :
%$$ \msb = \Upsilon_{s,\beta}^{1}(k_{1},k_{2}) \oplus \Upsilon_{s,\beta}^{2}(K_{1},K_{2}) \oplus \Upsilon_{s,\beta}^{3}(K_{1},K_{2})\,,$$
%where
%\begin{align*}
%\Upsilon_{s,\beta}^{1}(K_{1},K_{2})  = \big\{ M \in \msb, M_{[a]}^{[b]} = 0 \;\; &\mbox{if} \;\; \max(w_{a},w_{b}) \leq K_{1} \min(w_{a},w_{b}) \big\}\,,\\
%\Upsilon_{s,\beta}^{2}(K_{1},K_{2})  = \big\{ M \in \msb, M_{[a]}^{[b]} = 0 \;\; &\mbox{if} \;\; \max(w_{a},w_{b}) > K_{1} \min(w_{a},w_{b}) \\ &\mbox{or} \;\max(w_{a},w_{b})\leq K_{2} \big\}\,,\\
%\Upsilon_{s,\beta}^{3}(K_{1},K_{2})  = \big\{ M \in \msb, M_{[a]}^{[b]} = 0 \;\; &\mbox{if} \;\; \max(w_{a},w_{b}) > K_{1} \min(w_{a},w_{b}) \\ &\mbox{or} \;\max(w_{a},w_{b})> K_{2} \big\}\,,
%\end{align*}
%and we prove the desired estimates according to this decomposition. 

Since we estimate the operator norm of $B_{[a]}^{[b]}$, we need to rewrite the definition \eqref{eqdelort} in a operator way :  denoting by $D_{[a]}$ the diagonal (square) matrix with entries $\mu_{j}$, for $j \in [a]$ and $D'_{[a]}$ the diagonal (square) matrix with entries $k\cdot \om  +\eps\mu_j$, for $j \in [a]$, equation \eqref{eqdelort} reads
\be \label{eqdelort2} D'_{[a]} B_{[a]}^{[b]} -  B_{[a]}^{[b]} D_{[b]} =  i A_{[a]}^{[b]}\,.\ee
Then we distinguish 3 cases:

\noindent {\bf Case 1} : suppose  that $a,b$ satsify $$ \max(w_{a},w_{b}) > K_{1} \min(w_{a},w_{b})\,$$
take for instance $w_{a} > K_{1} w_{b}$. Then for $j\in[a]$
\be
| k\cdot \om \ +\eps\mu_j| \geq  w_{a} - \frac 14 - N|\omega| \geq \frac12 w_{a}\,,\ee
for 
\be\label{k1} K_{1}\geq4N|\omega| ,\ee
 that proves that $D'_{[a]}$ is invertible and gives an upper bound for the operator norm of its inverse. Then \eqref{eqdelort2} is equivalent to
\be B_{[a]}^{[b]} - {D'_{[a]}}^{-1} B_{[a]}^{[b]} D_{[b]} = i {D'_{[a]}}^{-1} A_{[a]}^{[b]}\,.\ee
Next consider the operator $\mathscr{L}^{1}_{[a]\times[b]}$ acting on matrices of size $[a]\times [b]$ such that
\be \mathscr{L}^{1}_{[a]\times[b]} \left( B_{[a]}^{[b]} \right) := {D'_{[a]}}^{-1} B_{[a]}^{[b]} D_{[b]}\,.\ee
We have 
\be \| \mathscr{L}^{1}_{[a]\times[b]} \left( B_{[a]}^{[b]} \right) \|  \leq  \frac{2 w_{b}}{w_{a }} \|B_{[a]}^{[b]}\| \leq  \frac{2}{K_{1}}  \|B_{[a]}^{[b]}\| \,,\ee
hence, in operator norm, $ \| \mathscr{L}^{1}_{[a]\times[b]} \| \leq \frac12$ if $K_{1}\geq 4$. Then the operator $\mathrm{Id} - \mathscr{L}^{1}_{[a]\times[b]}$ is invertible and
\begin{eqnarray*}
\| B_{[a]}^{[b]} \| &\leq& \| \left( \mathrm{Id} -  \mathscr{L}_{[a]\times[b]} \right)^{-1} \| \|  i {D'_{[a]}}^{-1} A_{[a]}^{[b]}\| \\
& \leq &  \frac{4}{w_{a}}  \|  A_{[a]}^{[b]}\| \,.
\end{eqnarray*}
But in case 1,  $1+|w_{a}-w_{b}| \leq 1+ w_{a}\leq  2 w_a$, therefore 
\be \label{delort1}\| B_{[a]}^{[b]} \| \leq 8 \frac{1}{1 + |w_{a}-w_{b}|}  \|  A_{[a]}^{[b]}\|.\ee

\noindent {\bf Case 2} :  suppose  that $a,b$ satisfy
$$ \max(w_{a},w_{b}) \leq K_{1} \min(w_{a},w_{b}) \; \mbox{and} \; \max(w_{a},w_{b})> K_{2} \,.$$
Notice that these two conditions imply that $$\min(w_a,w_b)\geq \frac{K_2}{K_1}.$$
We define the  square matrix $\tilde{D}_{[a]}=   w_{a} \mathbf{1}_{[a]}$, where $\mathbf{1}_{[a]}$ is the identity matrix. Then 
\be \| D_{[a]} - \tilde{D}_{[a]} \| \leq \frac {C_{\mu}}{w_{a}^{\delta}}\,,\ee
and equation \eqref{eqdelort} may be rewritten as
\be \label{amel2} \mathscr{L}^{2}_{[a]\times[b]} \left( B_{[a]}^{[b]} \right) -  (\tilde{D}_{[a]} -  D_{[a]})B_{[a]}^{[b]}  + B_{[a]}^{[b]} (\tilde{D}_{[b]} -  D_{[b]})= A_{[a]}^{[b]}  \,,\ee
where we denote by $\mathscr{L}^{2}_{[a]\times[b]}$  the operator acting on matrices of size $[a]\times [b]$ such that
\be \mathscr{L}^{2}_{[a]\times[b]} \left( B_{[a]}^{[b]} \right) :=  \left(  k\cdot\om +  w_{a} - w_{b}\right) B_{[a]}^{[b]} \,.
\ee
This dilation is invertible and \eqref{hypdelort} then gives, in operator norm,
\be \| \left(  \mathscr{L}^{2}_{[a]\times[b]} \right)^{-1}\| \leq  \frac{1}{\kappa (1 + |w_{a}-w_{b}|)}\,.\ee
This allows to write \eqref{amel2} as
\be \label{amel3} B_{[a]}^{[b]} - \left(  \mathscr{L}^{2}_{[a]\times[b]} \right)^{-1} \mathscr{K}_{[a]\times[b]}\left(B_{[a]}^{[b]} \right) =  \left(  \mathscr{L}^{2}_{[a]\times[b]} \right)^{-1}\left(A_{[a]}^{[b]}\right)\,,\ee
where $ \mathscr{K}_{[a]\times[b]}\left(B_{[a]}^{[b]} \right) =  (\tilde{D}_{[a]} -  D_{[a]})B_{[a]}^{[b]}  - B_{[a]}^{[b]} (\tilde{D}_{[b]} -  D_{[b]}) $. We have, thanks to \eqref{hypdelort0}, in operator norm,
\be \|\mathscr{K}_{[a]\times[b]} \| \leq C_{\mu}\left( \frac 1{w_{a}^{\delta}}+\frac 1{w_{b}^{\delta}}\right) \leq C_{\mu}\big(\frac{K_1}{K_{2}}\big)^{\delta}\,.\ee
Then for 
\be \label{k2}K_{2}\geq K_1(\frac{2C_{\mu}}{\kappa})^{1/\delta},\ee
 the operator $\mathrm{Id} - (\mathscr{L}^{2}_{[a]\times[b]})^{-1}\mathscr{K}_{[a]\times[b]}$ is invertible and from \eqref{amel3} we get 
\begin{eqnarray*} \| B_{[a]}^{[b]} \|  &=& \|\left( \mathrm{Id} - (\mathscr{L}^{2}_{[a]\times[b]})^{-1}\mathscr{K}_{[a]\times[b]}\right)^{-1} \| \| \left(  \mathscr{L}^{2}_{[a]\times[b]} \right)^{-1}\left(A_{[a]}^{[b]}\right)\| \\
& \leq & 2  \| \left(  \mathscr{L}^{2}_{[a]\times[b]} \right)^{-1}\left(A_{[a]}^{[b]}\right)\|\,,
\end{eqnarray*}
Hence in this case 
\be \label{delort2}\| B_{[a]}^{[b]} \| \leq \frac{2}{\kappa(1+|w_a-w_b|)}\| A_{[a]}^{[b]} \|\,.\ee

\noindent {\bf Case 3} : suppose  that $a,b\in\L$ satisfy
$$ \max(w_{a},w_{b}) \leq K_{1} \min(w_{a},w_{b}) \; \mbox{and} \; \max(w_{a},w_{b})\leq K_{2} \,.$$
In that case the size of the blocks are less than $K^d_2$ and we have
\be |B_{j}^{l}| = \left| \frac i{\langle k,\om(\r)\rangle \ +\eps\mu_j-\mu_l}\right| |{A}_j^l | \leq \frac{1}{\kappa(1+|w_{a}-w_{b}|) } |{A}_j^l |\ee
A majoration of the coefficients gives a poor majoration of the operator norm of a matrix, but it is sufficient here:
 \be \label{delort3} \| B_{[a]}^{[b]} \| \leq \frac{K_2^{d/2}}{\kappa(1+|w_{a}-w_{b}|) }\|A_{[a]}^{[b]} \|\,.
 \ee
Collecting \eqref{delort1}, \eqref{delort2} and \eqref{delort3} and taking into account \eqref{k1}, \eqref{k2} leads to the  result.
\endproof

\end{document}